\documentclass[MSNbibl,number,citesort,dvips]{arxbj}
\usepackage{mathrsfs,algorithm}
\usepackage{graphicx}

\aid{0}
\volume{19}
\issue{1}
\pubyear{2013}
\firstpage{64}
\lastpage{92}
\doi{10.3150/11-BEJ392}

\makeatletter
\def\d{\mathrm{d}}
\def\la{\langle}
\def\ra{\rangle}
\newtheorem{theorem}{Theorem}[section]
\newtheorem{proposition}[theorem]{Proposition}
\newremark{remark}[theorem]{Remark}
\newtheorem{lemma}[theorem]{Lemma}
\makeatother

\begin{document}
\begin{frontmatter}

\title{A Bayesian nonparametric approach to modeling market share dynamics}
\runtitle{Bayesian modeling of market dynamics}

\begin{aug}
%%%% inicialai - be tarpu
\author{\fnms{Igor} \snm{Pr\"unster}\thanksref{e1}\ead[label=e1,mark]{igor@econ.unito.it}} \and
\author{\fnms{Matteo} \snm{Ruggiero}\corref{}\thanksref{e2}\ead[label=e2,mark]{matteo.ruggiero@unito.it}}
\runauthor{I. Pr\"unster and M. Ruggiero} %% auto
\address{Dipartimento di Statistica e Matematica Applicata \& Collegio
Carlo Alberto,\\
Universit\`a degli Studi di Torino, C.so Unione Sovietica 218/bis,
10134 Torino, Italy.\\ \printead{e1,e2}}
\end{aug}

% HISTORY:
\received{\smonth{5} \syear{2010}}
\revised{\smonth{1} \syear{2011}}

% ABSTRACT
%
\begin{abstract}
We propose a flexible stochastic framework for modeling the market
share dynamics over time in a multiple markets setting, where firms
interact within and between markets. Firms undergo stochastic
idiosyncratic shocks, which contract their shares, and compete to
consolidate their position by acquiring new ones in both the market
where they operate and in new markets. The model parameters can
meaningfully account for phenomena such as barriers to entry and exit,
fixed and sunk costs, costs of expanding to new sectors with different
technologies and competitive advantage among firms. The construction is
obtained in a Bayesian framework by means of a collection of
nonparametric hierarchical mixtures, which induce the dependence
between markets and provide a generalization of the Blackwell--MacQueen
P\'olya urn scheme, which in turn is used to generate a partially
exchangeable dynamical particle system. A Markov Chain Monte Carlo
algorithm is provided for simulating trajectories of the system, by
means of which we perform a simulation study for transitions to
different economic regimes.
Moreover, it is shown that the infinite-dimensional properties of the
system, when appropriately transformed and rescaled, are those of a
collection of interacting Fleming--Viot diffusions.
\end{abstract}

% KEYWORDS
%

\end{frontmatter}

%s1 #&#
%s1 ###
\section{Introduction}

The idea of explaining firm dynamics by means of a stochastic model for
the market evolution has been present in the literature for a long
time. However, only recently, firm-specific stochastic elements have
been introduced to generate the dynamics. Jovanovic~\cite{J82} was the
first to formulate an equilibrium model where stochastic shocks are
drawn from a distribution with known variance and firm-specific mean,
thus determining selection of the most efficient. Later Ericson and
Pakes~\cite{EP85} provide a stochastic model for industry behavior
which allows for heterogeneity and idiosyncratic shocks, where firms
invest and the stochastic outcome determines the firm's success, thus
accounting for a selection process which can lead to the firm's exit
from the market.
Hopenhayn~\cite{H92} performs steady state analysis of a dynamic
stochastic model which allows for entry, exit and heterogeneity.
In~\cite{S07} a stochastic model for market share dynamics, based on
simple random walks, is introduced. The common feature of this
non-exhaustive list is that, despite the mentioned models being
inter-temporal and stochastic, the analysis and the explicit
description of the model dynamics are essentially done at equilibrium,
thus projecting the whole construction onto a static dimension and
accounting for time somehow implicitly.
%aggregate point of view, are complex dynamical systems, so their
%general analysis or even numerical computation turns out to be a
%complicated task''. It is the \emph{curse of dimensionality}, which
%prevents researchers from using realistic inter-temporal models which
%describe rich environments but are analytically intractable. Hence,}
Indeed, the researcher usually finds herself before the choice between
a dynamic model with a representative agent and a steady-state analysis
of an equilibrium model with heterogeneity. Furthermore, relevant for
our discussion are two technical difficulties with reference to
devising stochastic models for market share dynamics: the
interdependence of market shares, and the fact that the distribution of
the size of shocks to each firm's share is likely to depend on that
firm's current share. As stated in~\cite{S07}, these together imply
that an appropriate model might be one in which the distribution of
shocks to each firm's share is conditioned on the full vector of market
shares in the current period.

The urge to overcome these problems from an aggregate perspective,
while retaining the micro dynamics, has lead to a recent tendency of
borrowing ideas from statistical physics for modeling certain problems
in economics and finance. A particularly useful example of these tools
is given by interacting particle systems, which are
arbitrary-dimensional models describing the dynamic interaction of
several variables (or particles). These allow for heterogeneity and
idiosyncratic stochastic features but still permit a relatively easy
investigation of the aggregate system properties. In other words, the
macroscopic behavior of the system is derived from the microscopic
random interactions of the economic agents, and these techniques allow
us to keep track of the whole tree of outcomes in an inter-temporal framework.
A recent example of such an approach is given in~\cite{DPetal}, where
interacting particle systems are used to model the propagation of
financial distress in a network of firms. Another example is \cite
{R09}, which studies limit theorems for the process of empirical
measures of an economic model driven by a large system of agents that
interact locally by means of mechanisms similar to what, in population
genetics, are called mutation and recombination.

Here we propose a Bayesian nonparametric approach for modeling market
share dynamics by constructing a stochastic model with interacting
particles which allows us to overcome the above mentioned technical
difficulties.
In particular, a nonparametric approach allows us to avoid any
unnecessary assumption on the distributional form of the involved
quantities, while a Bayesian approach naturally incorporates
probabilistic clustering of objects and features conditional predictive
structures, easily admitting the representation of agents' interactions
based on the current individual status. Thus, with respect to the
literature on market share dynamics, we model time explicitly, instead
of analyzing the system at equilibrium, while retaining heterogeneity
and conditioning on the full vector of market shares. And despite the
different scope, with respect to the particles approach in~\cite{R09},
we instead consider many subsystems with interactions among each other
and thus obtain a vector of dependent continuous-time processes. In
constructing the model, the emphasis will be on generality and
flexibility, which necessarily implies a certain degree of stylization
of the dynamics. However, this allows the model to be easily adapted to
represent diverse applied frameworks, such as, for example, population
genetics, by appropriately specifying the corresponding relevant
parameters. As a matter of fact, we will follow the market share
motivation throughout the paper, with the parallel intent of favoring
intuition behind the stochastic mechanisms. A completely micro-founded
economic application will be provided in a follow-up paper \cite
{PRT11}. However, besides the construction, the present paper includes
an asymptotic distributional result which shows weak convergence of the
aggregate system to a collection of dependent diffusion processes. This
is a result of independent mathematical interest, relevant, in
particular, for the population genetics literature, where our
construction can be seen as a countable approximation of a system of
Fleming--Viot diffusions with mutation, selection and migration
(see~\cite{DG99}). Appendix A includes some basic material on
Fleming--Viot processes.

Finally, it is worth mentioning that our approach is also allied to
recent developments in the Bayesian nonparametric literature: although
structurally different, our model has a natural interpretation within
this field as belonging to the class of dependent processes, an
important line of research initiated in the seminal papers of \cite
{ME99,ME00}. Among others, we mention interesting dependent models
developed in~\cite{DIetal04,DGG07,DP08,PGG09,TMJ10} where one can find
applications to epidemiology, survival analysis and functional data
analysis. See the monograph~\cite{HHMW} for a recent review of the
discipline. Although powerful and flexible, Bayesian nonparametric
methods have not yet been extensively exploited for economic
applications. Among the contributions to date, we mention \cite
{MW05,LS2,GS10} for financial time series,~\cite{GS06,G07} for
volatility estimation,~\cite{LS} for option pricing and \cite
{BHH08,DJL10} for discrete choice models,~\cite{GS04} for stochastic
frontier models. With respect to this literature, the proposed
construction can be seen as a dynamic partially exchangeable array, so
that the dependence is meant both with respect to time and in terms of
a vector of random probability measures.

To be more specific, we introduce a flexible stochastic model for
describing the time dynamics of the market concentration in several
interacting, self-regulated markets. A~potentially infinite number of
companies operate in those markets where they have a positive share.
Firms can enter and exit a market, and expand or contract their share
in competition with other firms by means of endogenous stochastic
idiosyncratic shocks. The model parameters allow for barriers to entry
and exit, costs of expansion in new markets (e.g., technology
conversion costs), sunk costs and different mechanisms of competitive advantage.
The construction is achieved by first defining an appropriate
collection of dependent nonparametric hierarchical models and deriving
a related system of interacting generalized P\`olya urn schemes. This
underlying Bayesian framework is detailed in Section~\ref{sec:model}.
The collection of hierarchies induces the dependence between markets
and allows us to construct, in Section~\ref{sec:constr}, a dynamic
system, which is driven by means of Gibbs sampling techniques \cite
{GS90} and describes how companies interact among one another within
and between markets over time. These undergo stochastic idiosyncratic
shocks that lower their current share and compete to increment it. An
appropriate set of parameters regulates the mechanisms through which
firms acquire and lose shares and determines the competitive selection
in terms of relative strengths as functions of their current position
in the market and, possibly, the current market configuration as a
whole. For example, shocks can be set to be random in general but
deterministic when a firm crosses upwards some fixed threshold, meaning
that some antitrust authority has fixed an upper bound on the market
percentage which can be controlled by a single firm, which is thus
forced away from the dominant position. The competitive advantage
allows for a great degree of flexibility, involving a functional form
with very weak assumptions.
In Section~\ref{inf-prop} the dynamic system is then mapped into a
measure-valued process, which pools together the local information and
describes the evolution of the aggregate markets. The system is then
shown to converge in distribution, under certain conditions and after
appropriate rescaling, to a system of dependent diffusion processes,
each with values in the space of probability measures, known as
interacting Fleming--Viot diffusions.
In Section~\ref{sec:simulation} two algorithms which generate sample
paths of the system are presented, corresponding to competitive
advantage directly or implicitly modeled. A simulation study is then
performed to explore dynamically different economic scenarios with
several choices of the model parameters, investigating the effects of
changes in the market characteristics on the economic dynamics.
Particular attention is devoted to transitions of economic regimes as
dependent on specific features of the market, on regulations imposed by
the policy maker or on the interaction with other markets with
different structural properties.
Finally, Appendix~\ref{appmA} briefly recalls some background material
on Gibbs sampling, Fleming--Viot processes and interacting
Fleming--Viot processes, while all proofs are deferred to Appendix~\ref{appmB}.

%%%%%%%%%%%%%%%%%%%%%%%%%%%%%%%%%%%%%%%%%%%%%%%%%%%%%%%%%%%%%%%%
%%%%%%%%%%%%%%%%%%%%%%%%%%%%%%%%%%%%%%%%%%%%%%%%%%%%%%%%%%%%%%%%

%s2 #&#
%s2 ###
\section{The underlying framework}\label{sec:model}

In this section we define a collection of dependent nonparametric
hierarchical models, which will allow a dynamic representation of the
market's interaction.
%The Dirichlet process, introduced in~\cite{F73}, is a random
%probability measure which generalises the Dirichlet distribution, and
%can be defined as follows. Consider a measurable space $(\mathbb{X},
%space and $\mathscr{X}$ is the associated Borel sigma algebra, and let
%$\alpha$ be a non negative finite measure on $(\mathbb{X},
%$\alpha$, henceforth denoted $\mu\sim\mathscr{D}(\cdot|\alpha)$, if
%for every $k\ge1$ and every measurable partition $B_{1},\ldots,B_{k}$
%of $\mathbb{X}$, the vector $(\mu(B_{1}),\ldots,\mu(B_{k}))$ has
%Dirichlet distribution with parameters $(\alpha(B_{1}),\ldots,
%More constructive definitions of the Dirichlet process have been
%proposed in the literature.~\cite{S84}, for example, provided a series
%representation given by what is usually referred to as

%Another is due to~\cite{BM73} and is obtained via sequences of
%exchangeable observations. This will be useful for our purposes and is
%briefly outlined here.
Let $\alpha$ be a finite non-null measure on a complete and separable
space $\mathbb{X}$ endowed with its Borel sigma algebra $\mathscr
{X}$, and consider the P\`olya urn for a continuum of colors, which
represents a fundamental tool in many constructions of Bayesian
nonparametric models. This is such that $X_{1}\sim\alpha(\cdot
)/\alpha(\mathbb{X})$, and, for $n\ge2$,
%
%e1 #&#
\begin{equation}\label{eq: Polya-urn}
X_{n}|X_{1},\ldots,X_{n-1}\sim
\frac{\alpha(\cdot)+\sum_{i=1}^{n-1}\delta_{X_{k}}(\cdot)}{\alpha
(\mathbb{X})+n-1},
\end{equation}
where $\delta_{y}$ denotes a point mass at $y$. We will denote the
joint law of a sequence $(X_{1},\ldots,X_{n})$ from (\ref
{eq: Polya-urn}) with $\mathcal{M}_{n}^{\alpha}$, so that
%
%e2 #&#
\begin{equation}\label{P_n}
\mathcal{M}_{n}^{\alpha} =\frac{\alpha}{\alpha(\mathbb{X})}\prod
_{i=2}^{n}\frac{\alpha+\sum_{k< i}\delta_{X_{k}}}{\alpha(\mathbb{X})+i-1}.
\end{equation}
In~\cite{BM73} it is shown that this prediction scheme is closely
related to the Dirichlet process prior, introduced by~\cite{F73}. A
random probability measure $P$ on $(\mathbb{X},\mathscr{X})$ is said
to be a Dirichlet process with parameter measure $\alpha$, henceforth
denoted $P\sim\mathscr{D}(\cdot|\alpha)$, if for every $k\ge1$ and
every measurable partition $B_{1},\ldots,B_{k}$ of $\mathbb{X}$, the
vector $(P(B_{1}),\ldots,P(B_{k}))$ has Dirichlet distribution with
parameters $(\alpha(B_{1}),\ldots,\alpha(B_{k}))$.
%the conditional distribution of $X_{n}$ given $(X_{1},\ldots,X_{n-1})$
%converges almost surely to a discrete random probability measure $\mu$
%which is a Dirichlet process $\mathscr{D}(\cdot|\alpha)$, and, given $
%distribution $\mu$, that is they are exchangeable.
%An interpretation of the role of the Dirichlet process in this
%framework is provided by de Finetti's representation theorem,
%according to which $\mathscr{D}(\cdot|\alpha)$ is the unique
%probability measure which allows the representation, for every $n\ge
%1$,
%where $\mathscr{P}(\mathbb{X})$ denote the set of Borel probability
%measures on $\mathbb{X}$. Such $\mathscr{D}(\cdot|\alpha)\in
%and is also identified with the distribution of the almost sure weak
%limit of the empirical measure $n^{-1}\sum_{i\le n}x_{i}$. See for
%example~\cite{A85}, Section 1.2. Furthermore the left hand side of (

Among the various generalizations of the P\`olya urn scheme (\ref{eq:
Polya-urn}) present in the literature, a recent extension given in
\cite{RW09} will be particularly useful for our construction.
Consider, for every $n\ge1$, the joint distribution
%
%e3 #&#
\begin{equation}\label{joint}
q_n(\d x_1,\ldots,\d x_n)\propto p_n(\d x_1,\ldots,\d x_n) \prod
_{k=1}^{n}\beta_n(x_k),
\end{equation}
where $\beta_n$ is a given bounded measurable function on $\mathbb
{X}$. A representation for (\ref{joint}) can be provided in terms of a
Dirichlet process mixture model~\cite{L84}. In particular, it
can be easily seen that when $p_{n}\equiv\mathcal{M}_{n}^{\alpha}$
in (\ref{joint}), the predictive distribution for $X_{i}$, given
$\mathbf{x}_{(-i)}=(x_{1},\ldots,x_{i-1},x_{i+1},\ldots,x_{n})$, is
%
%For each $n$ even (which we assume henceforth), consider two vectors
%$(x_1,\ldots,x_n)$ and $(z_1,\ldots,z_n)$, whose coordinates are
%points in a complete separable metric space $\mathbb{X}$ and on the
%line respectively, endowed with their respective Borel sigma algebras.
%Denote with $P_n$ a pairing of $\{1,\ldots,n\}$ such that, given
%$P_n$, $k$ is paired with $j_k$. The distribution of the pairing $
%$P_n$ is $\tilde p_n(z_1,\ldots,z_n| x_1,\ldots,x_n,P_n)=\prod_{k} n
%K_n(z_k,z_{j_k}|x_k,x_{j_k})$, where the product is taken over $n/2$
%terms, hence covering all pairs. The density $K_n$ is assumed to be
%chosen for all $n$, and such that it equals $\beta_n(x_k,x_{j_k})$
%when computed in $(z_k= 1,z_{j_k}= 1)$, with $\beta_n\in B_{
%1|x_k,x_{j_k})=\beta_n(x_k,x_{j_k})$.
%Finally let $\{\mathscr{D}_{\alpha_{n}},n\in\mathbb{N}\}\subset
%processes, whose indexing measures $\alpha_{n}$ will be defined later.
%For each $n$, we take the vector $(x_1,\ldots,x_n)$ to be such that
%$x_i\sim\mu_{n}$ for every $i$, conditionally on $\mu_{n}\sim
%written, for every $n\ge1$, in the hierarchical form
% x_1,\ldots,x_n| \mu_{n}\stackrel{\mathrm{i.i.d.}}{\sim}& \mu_{n}
% \mu_{n}\sim& \mathscr{D}_{\alpha_{n}}\nonumber\\
% P_n\sim& \Pi_n\nonumber.
% \end{eqnarray}
%In particular, the above family of hierarchies defines a lower
%triangular matrix of $x$'s, given by $x_{1}|\mu_{1}\sim\mu_{1}$, then
%$x_{1},x_{2}|\mu_{2}\sim\mu_{2}$ and so on, where the $n$-th row can
%be thought of as the projections of an infinite exchangeable sequence
%conditionally from $\mu_{n}$, with $\mu_{n}\sim\mathscr{D}_{
%
%e4 #&#
\begin{equation}\label{pred}
q_{n,i}\bigl(\d x_{i}| \mathbf{x}_{(-i)}\bigr)
\propto\beta_n(x_{i})\alpha(\d x_{i})+\sum_{k\ne i}^n\beta
_n(x_{i})\delta_{x_k}(\d x_{i}).
\end{equation}
This can be thought of as a weighted version of (\ref{eq: Polya-urn}),
which is recovered when $\beta_{n}\equiv1$ for all $n\ge1$.
A more general version of (\ref{pred}) can be obtained by making
$\beta_{n}$ depend on the whole vector and thus allowing for a broad
range of interpretations. See the discussion following (\ref{pred2})
for this and for a more detailed interpretation of (\ref{pred}) in the
context of the present paper.

Consider now the following setting. For each $n$, let $(X_{1},\ldots
,X_{n})\in\mathbb{X}^{n}$ be an $n$-sized sample from $\mathcal
{M}_{n}^{\alpha}$, and let
%
%e5 #&#
\begin{equation}\label{eq: alpha_x}
\alpha_{x_{1},\ldots,x_{n}}(\d y)
=\alpha(\d y)+\sum_{k=1}^{n}\delta_{x_{k}}(\d y).
\end{equation}
%
%$\alpha_{x_{1},\ldots,x_{n}}:\mathbb{X}^{n}\times\mathscr{X}
%and
%x_{1},\ldots,\d x_{n})=\frac{\alpha(\d y)}{\alpha(\mathbb{X})},
% n\ge1.
Define the double hierarchy
%
%e6 #&#
\begin{eqnarray}\label{hierarchy3}
X_{1},\ldots,X_{n}| P&\stackrel{\mathrm{i.i.d.}}{\sim}& P, \qquad
P \sim\mathscr{D}(\cdot|\alpha),\nonumber
\\[-8pt]
\\[-8pt]
Y_{1},\ldots,Y_{n}| Q_{n}&\stackrel{\mathrm{i.i.d.}}{\sim}&Q_{n},
\qquad
Q_{n} \sim\mathscr{D}(\cdot|\alpha_{x_{1},\ldots,x_{n}}).
\nonumber
\end{eqnarray}
Here $(X_{1},\ldots,X_{n})$ are drawn from a Dirichlet process $P\sim
\mathscr{D}(\cdot|\alpha)$ and $(Y_{1},\ldots,Y_{n})$, given
$(X_{1},\ldots,X_{n})$, are drawn from a Dirichlet process
$Q_{n}:=Q|(X_{1},\ldots,X_{n})\sim\mathscr{D}(\cdot|\alpha
_{X_{1},\ldots,X_{n}})$. %, with %$\alpha_{x_{1},\ldots,x_{n}}$ as in (
%values of $(X_{1},\ldots,X_{n})$.
%Alternatively, $(Y_{1},\ldots,Y_{n})$ is drawn from the mixture of
%Dirichlet processes
%Q_{n}\sim\int_{\mathbb{X}^{n}}\mathscr{D}(\cdot|\alpha_{x_{1},
%as defined in~\cite{A74},
% Also take $\alpha$ to be
%x)p(y_{1},\ldots,y_{n})\d y_{1}\dots\d y_{n}
%where $p(y_{1},\ldots,y_{n})$ is the joint distribution of $(y_{1},
%Analogously to (\ref{definetti}),
It can be easily seen that
the joint law of $(Y_{1},\ldots,Y_{n})$ conditional on $(X_{1},\ldots
,X_{n})$ is $\mathcal{M}_{n}^{\alpha_{x_{1},\ldots,x_{n}}}$, with
$\mathcal{M}_{n}^{\alpha}$ as in (\ref{P_n}).
%See Proposition 1 and Corollary 3.2 and 3.2' in~\cite{A74}.
%%p_{n}(\d y_{1},\ldots,\d y_{n}| x_{1},\ldots,x_{n})
%Denote $\mathbf{x}=(x_{1},\ldots,x_{n})$ and $\mathbf{y}=(y_{1},
%expectation and Fubini's theorem from (\ref{hierarchy3}) we have that
%(with loose notation)
%p(\d\mathbf{x})=&\int_{\mathscr{P}(\mathbb{X})} p(\d\mathbf{x}|\mu)p(
%=&\int_{\mathscr{P}(\mathbb{X})} p(\d\mathbf{x}|\mu)\int_{\mathbb{X}}
%p(\d\mu|\mathbf{y})p(\d\mathbf{y})\\
%=&\int_{\mathbb{X}}\int_{\mathscr{P}(\mathbb{X})} p(\d\mathbf{x}|\mu)
%p(\d\mu|\mathbf{y})p(\d\mathbf{y})\\
%=& \mathbb{E}_{\mathbf{y}} (\int_{\mathscr{P}(\mathbb{X})} p(\d
%from which, conditioning on $(y_{1},\ldots,y_{n})$,
%p(\d\mathbf{x}|\mathbf{y})=&\int_{\mathscr{P}(\mathbb{X})} p(\d
%=&\int_{\mathscr{P}(\mathbb{X})} \mu^{n}(\d\mathbf{x})\mathscr{D}(\d
%which is analogous to (\ref{definetti}) with a different parameter.
%From this, (\ref{marg-x|y}) follows intuitively.
%
%$\alpha=\theta\nu_{0}$, where $\theta=\alpha(\mathbb{X})$ and $\nu_{0}=
The following result, stated here for ease of reference, can be found
in~\cite{WM03}.

%le2.1 #&#
\begin{lemma}\label{lemma: P_2n}
Let $\mathcal{M}_{n}^{\alpha}$ be as in (\ref{P_n}). Then
%
%e7 #&#
\begin{equation}\label{eq: marginal-requirement}
\int_{\mathbb{X}^{n}}\mathcal{M}_{n}^{\alpha_{x_{1},\ldots
,x_{n}}}(\d y_{1},\ldots,\d y_{n})\mathcal{M}_{n}^{\alpha}(\d
x_{1},\ldots,\d x_{n})=
\mathcal{M}_{n}^{\alpha}(\d y_{1},\ldots,\d y_{n}).
\end{equation}
\end{lemma}

In particular, Lemma~\ref{lemma: P_2n} yields a certain symmetry in
(\ref{hierarchy3}), so that we could also state that the joint law of
$(X_{1},\ldots,X_{n})$ conditional on $(Y_{1},\ldots,Y_{n})$ is
$\mathcal{M}_{n}^{\alpha_{y_{1},\ldots,y_{n}}}$.
Denote $\mathbf{x}=(x_{1},\ldots,x_{n})$ and extend (\ref{joint}) to
\begin{eqnarray*}%\label{}
q_n(\d\mathbf{x})\propto p_n(\d\mathbf{x}) \prod_{k=1}^{n}\beta_n(x_k),
\qquad
q_n(\d\mathbf{y})\propto p_n(\d\mathbf{y}) \prod_{k=1}^{n}\beta
_n(y_k).\nonumber
\end{eqnarray*}
From (\ref{pred}), when $(X_1,\ldots,X_n)$ and $(Y_1,\ldots,Y_n)$
come from (\ref{hierarchy3}) we have for $1\le i\le n$
%
%e8 #&#
\begin{equation}\label{pred-2col}
q_{2n,i}\bigl(\d x_{i}| \mathbf{x}_{(-i)},\mathbf{y}\bigr)
\propto
\beta_n(x_{i}) \alpha_{y_{1},\ldots,y_{n}}(\d x_{i})+\sum_{k\ne
i}^n\beta_n(x_k)\delta_{x_k}(\d x_{i})
\end{equation}
and similarly for $y_{i}$.
%This predictive distribution will be crucial for the specification of
%the model, constructed in Section~\ref{sec:constr}, which describes
%the market dynamics. %In particular, if we let each observation
%$x_{i}$ be an $n$-th fraction of the market identified by $x$, then (
%given the current configuration of two markets $x$ and $y$.
It is now straightforward to iterate the above argument and allow for
an arbitrary number of dependent hierarchies. Denote $\mathbf
{x}^{r}=(x^{r}_{1},\ldots,x^{r}_{n})$ and $\alpha_{\mathbf
{x}^{r}}=\alpha_{x^{r}_{1},\ldots,x^{r}_{n}}$, where $r,r',r''$
belong to some finite index set $\mathcal{I}$, whose cardinality is
denoted $\#\mathcal{I}$. Then, for every $n\ge1$, let
%
%e9 #&#
\begin{eqnarray}\label{hierarchy4}
\begin{array}{rcl@{ \qquad}rcl}
\mathbf{X}^{r}| P& \stackrel{\mathrm{i.i.d.}}{ \sim} &P,&
P& \sim&\mathscr{D}(\cdot|\alpha),\\\noalign{\vspace*{5pt}}
\mathbf{X}^{r}| P^{r}& \stackrel{\mathrm{i.i.d.}}{ \sim} &P^{r},&
P^{r}& \sim&\mathscr{D}(\cdot|\alpha_{\mathbf{x}^{r}}),\\\noalign{\vspace*{5pt}}
\mathbf{X}^{r''}| P^{r',r''}& \stackrel{\mathrm{i.i.d.}}{ \sim}
&P^{r',r''},&
P^{r',r''}& \sim&\mathscr{D}(\cdot|\alpha_{\mathbf{x}^{r},\mathbf
{x}^{r'}}),\\\noalign{\vspace*{5pt}}
&\vdots&&&\vdots&
\end{array}
\end{eqnarray}
where the dimension subscript $n$ has been suppressed in $\mathbf
{X}^{r},\mathbf{X}^{r'},\mathbf{X}^{r''},\ldots,$ for notational
simplicity. Denote now with
%
%e10 #&#
\begin{equation}\label{eq: system-dimension}
D_{n}=n\cdot\#\mathcal{I}
\end{equation}
the total number of components in (\ref{hierarchy4}). The joint law of
the $D_{n}$ items in (\ref{hierarchy4}) can be written
%
%e11 #&#
\begin{equation}\label{joint_particles}
%p_{_{D_{n}}}(\{\d\mathbf{x}^{r},r\in\mathcal{I}\})=
\mathcal{M}_{n}^{\alpha}(\d\mathbf{x}^{r})
\mathcal{M}_{n}^{\alpha_{\mathbf{x}^{r}}}(\d\mathbf{x}^{r'})
\mathcal{M}_{n}^{\alpha_{\mathbf{x}^{r},\mathbf{x}^{r'}}}(\d
\mathbf{x}^{r''})\cdots,
\end{equation}
where, in view of Lemma~\ref{lemma: P_2n}, (\ref{joint_particles}) is
invariant with respect to the order of $r,r',r'',\ldots.$
%$\prec$ is some ordering of $\mathcal{I}$ and
% \{
%
With a slight abuse of notation, define
%
%e12 #&#
\begin{eqnarray}\label{eq: Xi-xi}
\mathcal{I}(-r)= \{\mathbf{x}^{r'}\dvt r'\in\mathcal{I}, r'\ne r\},
\end{eqnarray}
to be the set of all system components without the vector $\mathbf
{x}^{r}$, and
%
%e13 #&#
\begin{eqnarray}\label{eq: Xi-x_i}
\mathcal{I}(-x^{r}_{i})=\{\mathbf{x}^{r'}\dvt r'\in\mathcal{I}\}
\setminus\{x^{r}_{i}\}
\end{eqnarray}
to be the set of all system components without the item $x_{i}^{r}$.
%With some abuse of notation, the symbol $\mathcal{I}(-r)$ will be also
%used to denote $\mathcal{I}\setminus\{r\}$.
Analogously to (\ref{pred-2col}) in this enlarged framework,
the predictive law for $x^{r}_{i}$, conditional on the rest of the
system, can be written
%
%e14 #&#
\begin{equation}\label{glob-pred2}
q_{_{D_{n},i}}(\d x^{r}_{i}| \mathcal{I}(-x^{r}_{i}))
\propto\beta_n(x^{r}_i) \alpha_{_{\mathcal{I}(-r)}}(\d
x^{r}_{i})+\sum_{k\ne i}^n\beta_n(x^{r}_k) \delta_{x^{r}_k}(\d x^{r}_{i}),
\end{equation}
%
%Here $\alpha_{_{\mathcal{I}(-r)}}:\mathbb{X}^{D_{n}-1}\times\mathscr{B(
where the interpretation of $\alpha_{_{\mathcal{I}(-r)}}$ is clear
from (\ref{eq: alpha_x}) and (\ref{eq: Xi-xi}).
Note that this predictive law reduces to (\ref{eq: Polya-urn}) when
$\beta_{n}\equiv1$ and $\alpha_{_{\mathcal{I}(-r)}}\equiv\alpha$.
%marginal-requirement}) becomes
Expression (\ref{glob-pred2}) will be the key for the definition of
the market dynamics by means of an interacting system of particles.
A~detailed interpretation for $q_{D_{n},i}$ will be provided in the
following section. See (\ref{pred2}) and the following discussion.

To conclude the section, it is worth noting that (\ref{hierarchy4})
generates a partially exchangeable array, where partial exchangeability
is intended in the sense of de Finetti (see, e.g.,~\cite{CR96}). That
is, if $r,r',r''$ identify rows, then the system components are
row-wise exchangeable but not exchangeable.

%%%%%%%%%%%%%%%%%%%%%%%%%%%%%%%%%%%%%%%%%%%%%%%%%%%%%%%%%%%%%%%%%
%%%%%%%%%%%%%%%%%%%%%%%%%%%%%%%%%%%%%%%%%%%%%%%%%%%%%%%%%%%%%%%%%

%s3 #&#
%s3 ###
\section{Dynamic models for market evolution}\label{sec:constr}

In this section we define a dynamical model for the temporal evolution
of the firms' market shares in multiple interacting markets. The model
can be regarded as a random element whose realizations are
right-continuous functions from $[0,\infty)$ to the space $(\mathbb
{X}^{D_{n}},\mathscr{X}^{D_{n}})$, $D_{n}\in\mathbb{N}$ being (\ref
{eq: system-dimension}), and we refer to it as a \emph{particle
system}, since it explicitly models the evolution of the share units,
or particles, in several markets. For ease of presentation, we approach
the construction by first considering a single market for a fixed
number $n$ of share units, and then extend it to a collection of
markets. The investigation of the asymptotic properties as
$n\rightarrow\infty$ is, instead, the object of Section~\ref{inf-prop}.

For any fixed $n\ge1$, consider a vector $\mathbf{x}=\mathbf
{x}^{(n)}=(x_{1},\ldots,x_{n})\in\mathbb{X}^{n}$, and let
$(x_{1}^{*},\ldots,x_{K_{n}}^{*})$ denote the $K_{n}\le n$ distinct
values in $\mathbf{x}$, with $x_{j}^{*}$ having multiplicity $n_{j}$.
The elements of $(x_{1}^{*},\ldots,x_{K_{n}}^{*})$ represent the
$K_{n}$ firms operating in the market at a given time. Here $x_{j}^{*}$
is a random label to be seen as a unique firm identifier.
The vector $\mathbf{x}$ represents the current market configuration,
carrying implicitly the information on the shares.
Namely, the fraction of elements in $\mathbf{x}$ equal to $x_{j}^{*}$
is the market share possessed by firm $j$. Here $n$ represents the
level of share fractionalization in the market. Dividing the market
into $n$ fractions is not restrictive, since any share can be
approximated by means of a sufficiently large $n$. See Remark \ref
{remark:fractionalization} below for a discussion of the implications
of this assumption on the computational costs.

Define now a Markov chain taking values in $\mathbb{X}^{n}$ as follows.
At each step an index $i$ is chosen from $\{1,\ldots,n\}$ with
probability $\gamma_{n,i}\ge0$ for $i=1,\ldots,n$, with $\sum
_{i=1}^{n}\gamma_{n,i}= 1$.
Equivalently, let $\gamma_{j}(\mathbf{n}_{n})$ be the probability
that firm $x_{j}^{*}$ loses an $n$th fraction of its market share at a
certain transition, where $\gamma_{j}(\mathbf{n}_{n})$ depends on the
frequencies $\mathbf{n}_{n}=(n_{1},\ldots,n_{K_{n}})$. That is, firm
$x_{j}^{*}$ undergoes a shock whose probability is idiosyncratic,
depending on the firm itself and on the current market configuration,
summarized by the vector of frequencies. Different choices of $\gamma
_{j}(\mathbf{n}_{n})$ reflect different market regulations, possibly
imposed by the policy maker.
We provide some examples:
\begin{enumerate}[(3)]
\item[(1)] $\gamma_{j}(\mathbf{n}_{n})=1/K_{n}$: neutrality. All
firms have equal probability of undergoing a shock;

\item[(2)] $\gamma_{j}(\mathbf{n}_{n})=n_{j}/n$:
firms with higher shares
are the weakest, with a flattening effect on the share distribution.
This parametrization is also useful in population genetics contexts,
where particles represent individuals;

\item[(3)] $\gamma_{j}(\mathbf{n}_{n})=(1-n_{j}/n)/(K_{n}-1)$ when
$K_{n}\ge2$:
firms with higher shares
are the strongest.
The probability of losing shares is decreasing in the firms' positions
in the market;

\item[(4)] $\gamma_{j}(\mathbf{n}_{n})=\mathbf{1}(\max_i n_{i}\le
nC)\tilde\gamma_{j}(\mathbf{n}_{n})+\mathbf{1}(n_{j}> nC)$ for some
constant $0<C<1$, where $\mathbf{1}(A)$ is the indicator function of
the event $A$. The probability of selecting $x_{j}^{*}$ is $\tilde
\gamma_{j}(\mathbf{n}_{n})$ provided
no firm controls more than %$x_{j}^{*}$ controls at most
$C\%$ of the market. If firm $x_{j}^{*}$ controls more than $C\%$
of the market, at the following step, $x_{j}^{*}$ is selected with
probability one. Thus $C$ is an upper bound imposed by the policy maker
to avoid dominant positions. Incidentally, there is a subtler aspect of
this mechanism which is worth commenting upon. It will be seen later
that there is positive probability that the same firm acquires the
vacant share again, but this only results in picking again $x_{j}^{*}$
with probability one, until the threshold $C$ is crossed downwards.
This seemingly anomalous effect can be thought of as the viscosity with
which a firm in a dominant position gets back to a legitimate status
when condemned by the antitrust authority, which in no real world
occurs instantaneously.
\end{enumerate}

Suppose now $x_{i}=x_{j}^{*}$ has been chosen in $\mathbf{x}$.
Once firm $x_{j}^{*}$ looses a fraction of its share, the next state of
the chain is obtained by sampling a new value for $X_{i}$ from (\ref
{pred}), leaving all other components unchanged. Hence the $i$th
fraction of share is reallocated, according to the predictive
distribution of $X_{i}|\mathbf{x}_{(-i)}$, either to an existing firm
or to a new one entering the market.

%re3.1 #&#
\begin{remark}\label{remark:gibbs sampler}
The above Markov chain can also be thought of as generated by a Gibbs
sampler on $q_{n}(\d x_{1},\ldots,\d x_{n})$. This consists of
sequentially updating one randomly selected component at a time in
$(x_{1},\ldots,x_{n})$ according to the component-specific \emph{full
conditional} distribution $q_{n,i}(\d x_{i}|\mathbf{x}_{(-i)})$. The
Gibbs sampler is a special case of a Metropolis--Hastings Markov chain
Monte Carlo algorithm, and, under some assumptions satisfied within the
above framework, yields a chain which is reversible with respect to
$q_{n}(\d x_{1},\ldots,\d x_{n})$, hence also stationary. See \cite
{GS90} for details and Appendix~\ref{appmA} for a brief account.
\end{remark}

Consider now an arbitrary collection of markets, indexed by
$r,r',r'',\ldots\in\mathcal{I}$, so that the total size of the
system is (\ref{eq: system-dimension}), and extend the construction as follows.
%At rate $\lambda_{n}$ a transition occurs.
At each transition, a market $r$ is selected at random with probability
$\varrho_{r}$, and a component of $(x^{r}_{1},\ldots,x^{r}_{n})$ is
selected at random with probability $\gamma^{r}_{n,i}$. The next state
is obtained by setting all components of the system, different from
$x^{r}_{i}$, equal to their previous state, and by sampling a new value
for $x^{r}_{i}$ from (\ref{glob-pred2}). Choose now
%
%e15 #&#
\begin{eqnarray}\label{alpha}
\alpha_{_{\mathcal{I}(-r)}}(\d y)= \theta\pi\nu_{0}(\d y)+\theta
(1-\pi)\sum_{r'\in\mathcal{I}}m(r,r')\mu_{r'}(\d y),
\end{eqnarray}
where $\theta>0$, $\pi\in[0,1]$, $\nu_{0}$ is a non-atomic
probability measure on $\mathbb{X}$,
%
%e16 #&#
\begin{equation}\label{xi-empirical}
\mu_{r'}=n^{-1}\sum_{i=1}^{n}\delta_{x^{r'}_{i}}
\end{equation}
and $m(r,r')\dvtx \mathcal{I}\times\mathcal{I}\rightarrow[0,1]$ is such that
%
%e17 #&#
\begin{equation}\label{m}
m(r,r)=0,
\qquad
\sum_{r'\in\mathcal{I}}m(r,r')=1.
\end{equation}
In this case (\ref{glob-pred2}) becomes
%
%e18 #&#
\begin{eqnarray}\label{pred2}
q_{_{D_{n},i}}(\d x^{r}_i| \mathcal{I}(-x^{r}_{i}))
&\propto& \theta\pi\beta_n(x^{r}_i)\nu_{0}(\d
x^{r}_{i})\nonumber
\\[-8pt]
\\[-8pt]
&& {} +\theta(1-\pi)\beta_n(x^{r}_i)\sum_{r'\in\mathcal
{I}}m(r,r')\mu_{r'}(\d x^{r}_{i})
+\sum_{k\ne i}^{n}\beta_n(x^{r}_k) \delta_{x^{r}_k}(\d x^{r}_i)
\nonumber
\end{eqnarray}
with normalizing constant
%x^{r}_{i})+\theta(1-\pi)\int\beta_n(x^{r}_i)\sum_{r'\in
%+\sum_{k\ne i}^n\beta_n(x^{r}_k).
$\bar q_{_{D_{n},i}}=\mathrm{O}(n)$ when $\beta_{n}=1+\mathrm{O}(n^{-1})$.
By inspection of (\ref{pred2}), there are three possible destinations
for the allocation of the vacant share:
%. Recall that $x^{r}_{i}$ denotes the $i$-th component of market $r$,
%$x^{r}_{j^{*}}$ the $j$-th distinct value, or cluster, in $(x^{r}_{1},
%the associated absolute frequency. Then a sample of size one from (
%
\begin{enumerate}[(iii)]
\item[(i)] A new firm is created and enters the market. The new
value of the location $x^{r}_{i}$ is sampled from $\nu_{0}$, which is
non-atomic, so that $x^{r}_{i}$ has (almost surely) never been
observed. Here $\nu_{0}$ is common to all markets. The possibility of
choosing different $\nu_{0,r}$, $r\in\mathcal{I}$, is discussed in
Section~\ref{sec:simulation} below.

\item[(ii)] A firm operating in the same market $r$ expands its
share. The location is sampled from the last term, which is a weighted
empirical measure of the share distribution in market $r$, obtained by
ignoring the vacant share unit $x^{r}_{i}$.

\item[(iii)] A firm operating in another market $r'$ either enters
market $r$ or expands its current position in $r$. The location is
sampled from the second term. In this case, an index $r'\neq r$ is
chosen according to the weights $m(r,\cdot)$; then within $r'$ a firm
$x^{r'}_{j^{*}}$ is chosen according to the weighted empirical measure
\[
\mu_{r'}(\d y)=n^{-1}\sum_{k=1}^{n}\beta_n(x^{r'}_k)\delta
_{x^{r'}_{k}}(\d y).
\]
If the cluster associated to $x^{r}_{i}$ has null frequency in the
current state, we have an entrance from $r'$; otherwise, we have a
consolidation in $r$ of a firm that operates, at least, on both those markets.
\end{enumerate}
We can now provide interpretation for the model parameters:
\begin{enumerate}[(b)]
\item[(a)] $\theta$ governs \emph{barriers to entry}: the lower the
$\theta$, the higher the barriers to entry, both for entrance of new
firms and for those operating in other markets.

\item[(b)] $\pi$ regulates \emph{sunk costs}: given $\theta$, a low
$\pi$ makes expansions from other sectors more likely than start-up of
new firms, and vice versa.

\item[(c)] $m(r,\cdot)$ allows us to set the costs of expanding to
different sectors. For example, it might represent \emph{costs of
technology conversion} a firm needs to sustain or some \emph
{regulation} constraining its ability to operate in a certain market.
Tuning $m(\cdot,\cdot)$ on the base of some notion of distance
between markets allows us to model these costs, so that a low $m(r,r')$
implies, say, that $r$ and $r'$ require very different technologies,
and vice versa.

\item[(d)] $\beta_{n}$ is probably the most flexible parameter of the
model, which, due to the minimal assumptions on its functional form
(see Section~\ref{sec:model}), can reflect different features of the
market, implying several possible interpretations. For example, it
might represent \emph{competitive advantage}.
%role of $\beta_{n}$ in (\ref{pred2}), with respect to the first and
%second terms, is negligible for large $n$.}
Since $\beta_{n}$ assigns different weights to different locations of
$\mathbb{X}$, the higher $\beta_{n}(x^{r}_{j^{*}})$, the more favored
is $x^{r}_{j^{*}}$ when competing with the other firms in the same
market. Here, and later, $x^{r}_{j^{*}}$ denotes the $j$th firm in
market $r$. It is, however, to be noted that setting $\beta\equiv1$
does not imply competitive neutrality among firms, as the empirical
measure implicitly favors those with higher shares.
More generally, observe that the model allows us to consider a weight
function of type $\beta_{n}(x^{r}_{k},\mu_{r})$, where $\mu_{r}$ is
the empirical measure of market $r$, making $\beta_{n}$ depend on the
whole current market configuration and on $x^{r}_{k}$ explicitly. This
indeed allows for multiple interpretations and to arbitrarily set how
firms relate to one another when competing in the same market. For
example, this more general parametrization allows us to model
neutrality among firms
by setting $\beta_{n}(x^{r}_{k},\mu_{r})=1/n_{j}^{r}$, with
$n_{j}^{r}$ being the number of share units possessed by firm $j$ in
market $r$.

\item[(e)] Weights $\gamma_{n,i}$ can model \emph{barriers to exit},
if appropriately tuned (see also points (1) to (4) above). For example,
setting $\gamma_{j}(\mathbf{n}_{n})$ very low (null) whenever
$n_{j}$, or $n_{j}/n$, is lower than a given threshold makes the exit
of firm $x^{r}_{j^{*}}$ very unlikely (impossible).
\end{enumerate}

The function $\beta_{n}$, in point (d) above, will represent the
crucial quantity which will be used for introducing explicitly the
micro-foundation of the model. However, we do not pursue this here
since we focus on generality and adaptability of the model. The
micro-foundation will be the object of a subsequent work.

%s4 #&#
%s4 ###
\section{Infinite dimensional properties}\label{inf-prop}

From a qualitative point of view the outlined discrete-time
construction would be enough for many applications. Indeed Section \ref
{sec:simulation} below presents two algorithms which generate
realizations of the system and are used to perform a simulation study,
based on the above description.
It is, however, convenient to embed the chain in continuous time, which
makes the investigation of its properties somewhat simpler
and leads to a result of independent mathematical interest. This will
enable us to show that an appropriate transformation of the continuous
time chain converges in distribution to a well-known class of processes
which possess nice sample path properties. To this end, superimpose the
chain to a Poisson point process with intensity $\lambda_{n}$, which
governs the waiting times between points of discontinuity.
The following proposition identifies the generator of the resulting
process under some specific assumptions which will be useful later.
Recall that the infinitesimal generator of a stochastic process $\{
Z(t),t\ge0\}$ on a Banach space $L$ is the linear operator $A$ defined by
\[
Af=\lim_{t\downarrow0}\frac{1}{t} \bigl[\mathbb{E}[f(Z(t))|Z(0)]-f(Z(0)) \bigr]
\]
with domain given by the subspace of all $f\in L$, for which the limit
exists. In particular, the infinitesimal generator carries all the
essential information about the process, since it determines the
finite-dimensional distributions.
Before stating the result, we need to introduce some notation. Let
$B(\mathbb{X})$ be the space of bounded measurable functions on
$\mathbb{X}$, and
$(\varrho_{r})_{r\in\mathcal{I}}$ be a sequence with values in the
corresponding simplex
%
%e19 #&#
\begin{equation}\label{simplex}
\Delta_{\#\mathcal{I}}= \biggl\{(\varrho_{r})_{r\in\mathcal{I}}\dvt \rho
_{r}\ge0, \forall r\in\mathcal{I}, \sum_{r\in\mathcal{I}}\varrho
_{r}=1 \biggr\}.
\end{equation}
Furthermore, let $q_{_{D_{n},i}}$ be as in (\ref{pred2}), with $\bar
q_{_{D_{n},i}}$ its normalizing constant, %$\alpha_{_{
and let
%
%e20 #&#
%e21 #&#
\begin{eqnarray}\label{sel_funct}
\beta_n(z)&=& 1+\sigma(z)/n,\qquad\sigma\in B(\mathbb{X}),\\\label
{constant}
C_{n,r,i}&=& \lambda_{n} \varrho_{r} \gamma_{n,i}^{r}/\bar q_{_{D_{n},i}}.
\end{eqnarray}
Define also the operators
%
%e22 #&#
%e23 #&#
%e24 #&#
\begin{eqnarray}\label{eta}
\eta_i(\mathbf{x}|z)&=&(x_1,\ldots,x_{i-1},z,x_{i+1},\ldots,x_n),\\
\label{mut-op}
M^{n}g(w)&=&\int[g(y)-g(w) ] \bigl(1+\sigma(y)/n \bigr)\nu_{0}(\d y),
\qquad g\in B(\mathbb{X}),\\\label{mig-op}
G^{n,r'}g(w)&=&\int[g(y)-g(w) ] \bigl(1+\sigma(y)/n \bigr)\mu_{r'}(\d y),
\qquad g\in B(\mathbb{X}),
\end{eqnarray}
and denote by
%
%e25 #&#
\begin{equation}\label{xi-oper}
\eta_{r_{i}},\qquad
M_{r_{i}}^{n}f,\qquad
G_{r_{i}}^{n,r'}f
\end{equation}
such operators as applied to the $i$th coordinate of those in ${\mathbf
{x}}$ which belong to $r$.
For instance, if $\mathbf
{y}=(y_{1}^{r'},y_{2}^{r},y_{3}^{r},y_{4}^{r'})$, where $y_{2},y_{3}$
belong to market $r$ and the others to $r'$,
%meaning that $(y_{3},y_{5},y_{6})$ belong to market $r$ and
%$(y_{1},y_{2},y_{4})$ to $r'$. T
then $\eta_{r_{2}}(\mathbf{y}|z)=\eta_{3}(\mathbf{y}|z)=
(y_{1}^{r'},y_{2}^{r},z,y_{4}^{r'})$.
%That is, for an operator $H$,
%H_{r_{i}}g(\mathbf{y})
%denote the operator $H$ applied to $g$ only as a function of the
%$i$-th of those components of $\mathbf{y}$ that belong to $r$.

%pr4.1 #&#
\begin{proposition}\label{theorem:particle-gen}
Let $X^{(D_{n})}(\cdot)=\{X^{(D_{n})}(t),t\ge0\}$ be the
right-continuous process with values in $\mathbb{X}^{D_{n}}$ which
updates
one component according to (\ref{pred2}) at each point of a Poisson
point process with intensity $\lambda_{n}$.
%
%one component at each point of discontinuity according to
%(\ref{pred2}).
Then $X^{(D_{n})}(\cdot)$ has infinitesimal generator, for
$f\in B(\mathbb{X}^{D_{n}})$, given by
%
%e26 #&#
\begin{eqnarray}\label{genallcol2}
A_{D_{n}}f(\mathbf{x})
&=&\sum_{r\in\mathcal{I}} \Biggl\{\theta\pi\sum_{i=1}^{n}C_{n,r,i}
M_{r_{i}}^{n}f(\mathbf{x})\nonumber\\
&&\hphantom{\sum_{r\in\mathcal{I}}
\Biggl\{}{}+\theta(1-\pi)\sum_{r'}m(r,r')\sum_{i=1}^{n}C_{n,r,i}
G_{r_{i}}^{n,r'}f(\mathbf{x})\nonumber
\\[-8pt]
\\[-8pt]
&&\hphantom{\sum_{r\in\mathcal{I}} \Biggl\{}{}
+\sum_{1\le k\ne i\le n}C_{n,r,i} [f(\eta_{r_{i}}(\mathbf
{x}|x^{r}_k))-f(\mathbf{x}) ]\nonumber\\
&&\hphantom{\sum_{r\in\mathcal{I}}
\Biggl\{}{}+\frac{1}{n}\sum_{1\le k\ne i\le n}C_{n,r,i} \sigma(x^{r}_k)
[f(\eta_{r_{i}}(\mathbf{x}|x^{r}_k))-f(\mathbf{x}) ] \Biggr\}.\nonumber
\end{eqnarray}
\end{proposition}

With respect to the market dynamics, generator (\ref{genallcol2}) can
be interpreted as follows. The first term governs the creation of new
firms, obtained by means of operator (\ref{mut-op}) which updates with
new values from $\nu_{0}$. The second regulates the entrance of firms
from other markets, via operator (\ref{mig-op}) and according to the
``distance'' kernel $m(\cdot,\cdot)$. The last two terms deal with
the expansion of firms in the same market. These parallel,
respectively, points (i), (iii) and (ii) above.

Consider now the probability-measure-valued system associated with
$X^{(D_{n})}(\cdot)$, that is, $Y^{(n)}(\cdot)=\{Y^{(n)}(t),t\ge0\}
$, where
%
%e27 #&#
\begin{equation}\label{Yn}
Y^{(n)}(t)=(\mu_{r}(t),\mu_{r'}(t),\ldots)
\end{equation}
and $\mu_{r}$ is as in (\ref{xi-empirical}). $Y^{(n)}(t)$ is thus the
collection of the empirical measures associated to each market, which
provides aggregate information on the share distributions at time $t$.
The following result identifies the generator of $Y^{(n)}(\cdot)$, for
which we need some additional notation.
Let
%
%e28 #&#
\begin{equation}\label{desc-factorial}
n_{[k]}=n(n-1)\cdots(n-k+1), \qquad
n_{[0]}=1.
\end{equation}
For every sequence $(r_{1},\ldots,r_{m})\in\mathcal{I}^{m}$, $m\in
\mathbb{N}$, and given $r\in\mathcal{I}$, define $k_{r}=\sum
_{j=1}^{m}\mathbf{1}(r_{j}=r)$ to be the number of elements in
$(r_{1},\ldots,r_{m})$ equal to $r$. Define also $\mu_{r}^{(k_{r})}$
and $\mu^{(m)}$ to be the probability measures
%
%e29 #&#
%e30 #&#
\begin{eqnarray}\label{mu-(k-xi)}
\mu_{r}^{(k_{r})}
&=& \frac{1}{n_{[k_{r}]}}\sum_{1\le i_{r,1}\ne\cdots\ne
i_{r,k_{r}}\le n}\delta_{(x^{r}_{i_{r,1}},\ldots
,x^{r}_{i_{r,k_{r}}})},\\
\label{mu-(m)}
\mu^{(m)}
&=& \prod_{r\in\mathcal{I}}\mu_{r}^{(k_{r})},
\end{eqnarray}
and let
%
%e31 #&#
\begin{equation}\label{phi}
\phi_{m}(\mu)=\int f\,\d\mu^{(m)},
\qquad
f\in B(\mathbb{X}^{m}).
\end{equation}
Finally, denote $\sigma_{r_{k}}(\cdot)=\sigma(x_{k}^{r})$ and
%
%e32 #&#
\begin{equation}\label{xi-oper2}
r_{k,i}=(r_{k},r_{i}),
\end{equation}
with $r_{i}$ as in (\ref{xi-oper}), and
define the map $\Phi_{ki}\dvtx B(\mathbb{X}^{n})\rightarrow B(\mathbb
{X}^{n-1})$ by
%
%e33 #&#
\begin{equation}\label{Phi}
\Phi_{ki}f(x_{1},\ldots,x_{n})=f(x_{1},\ldots
,x_{i-1},x_{k},x_{i+1},\ldots,x_{n}).
\end{equation}

%pr4.2 #&#
\begin{proposition}\label{theorem:mv-generator}
Let $Y^{(n)}(\cdot)$ be as in (\ref{Yn}).
%Then, for $f\in B(\mathbb
%{X}^{m})$, $m\le D_{n}$ and under the hypothesis and notation of
%Proposition~\ref{theorem:particle-gen}, the generator of
%$Y^{(n)}(\cdot)$ is, for any $m\le D_{n}$,
%
Then, for $\phi_{m}(\mu)$ as in (\ref{phi}), $m\leq D_n$ and under the
hypothesis and notation of Proposition~\ref{theorem:particle-gen},
the generator of $Y^{(n)}(\cdot)$ is
%
%e34 #&#
\begin{eqnarray}\label{mv-gen-(m)}
\mathbb{A}_{D_{n}}\phi_{m}(\mu)
&=&\sum_{r\in\mathcal{I}} \Biggl\{\theta\pi\sum_{i=1}^mC_{n,r,i}\int
M_{r_{i}}^{n}f\,\d\mu^{(m)}\nonumber\\
%&+\theta\pi(n-k_{r})\sum_{i=1}^mC_{n,r,i}\la M_{r_{i}}^{n,m+1}f,
&&\hphantom{\sum_{r\in\mathcal{I}} \Biggl\{}{}+\theta(1-\pi)\sum_{r'}m(r,r')\sum_{i=1}^mC_{n,r,i} \int
G_{r_{i}}^{n,r'}f\,\d\mu^{(m)}\nonumber\\
% &+\frac{\theta(1-\pi)(n-k_{r})}{n-1}\sum_{r'}a(r,r')\sum_{i=1}^m
%C_{n,r,i}\la G_{r_{i}}^{n,m+1,r'}f,\mu^{(m+1)}\ra\nonumber\\
&&\hphantom{\sum_{r\in\mathcal{I}} \Biggl\{}{}+\sum_{1\le k\ne i\le m}C_{n,r,i}\int(\Phi_{r_{k,i}}f- f)\,\d\mu
^{(m)}\\
&&\hphantom{\sum_{r\in\mathcal{I}} \Biggl\{}{}+\frac{1}{n}\sum_{i=1}^m\sum_{k\ne i}^{k_{r}}C_{n,r,i}\int\sigma
_{r_{k}}(\cdot)(\Phi_{r_{k,i}}f-f)\,\d\mu^{(m)}\nonumber\\
&&\hphantom{\sum_{r\in\mathcal{I}} \Biggl\{}{}+\frac{n-k_{r}}{n}\sum_{i=1}^mC_{n,r,i} \int
\sigma_{m+1}(\cdot)(\Phi_{m+1,r_{i}}f-f)\,\d\mu^{(m+1)}
%&+\frac{n-k_{r}}{n}\sum_{i=1}^m\sum_{k\ne i}^{k_{r}}C_{n,r,i}\la
%&+\frac{2(n-k_{r})(n-k_{r}-1)}{n(n-1)}\sum_{i=1}^mC_{n,r,i}\la
\Biggr\}.\nonumber
\end{eqnarray}
%
%with domain $\mathscr{D}(\mathbb{A}_{D_{n}})=\{\phi_{m}(\mu):\phi_{m}(
\end{proposition}

%This is the generator of the $(\mathscr{P}(\mathbb{X}))^{\#
%shares, when each market is divided in $n$ fractions.

%This form of the measure-valued system generator will be needed in
%Section~\ref{sec:limit} when the asymptotic properties of the system
%will be investigated for $n$ tending to infinity. Section
%processes, which is a collection of interacting diffusions. Then in
%Section~\ref{sec:limit} we show that the system of market share
%distributions converges weakly to a system of interacting
%Fleming--Viot processes.

The interpretation of (\ref{mv-gen-(m)}) is similar to that of (\ref
{genallcol2}), except that (\ref{mv-gen-(m)}) operates on the product
space $\mathscr{P}(\mathbb{X})^{\#\mathcal{I}}$ instead of the
product space of particles.
Let $\mathscr{P}^{n}(\mathbb{X})\subset\mathscr{P}(\mathbb{X})$ be
the set of purely atomic probability measures on $\mathbb{X}$ with
atom masses proportional to $n^{-1}$, $D_{\mathscr{P}(\mathbb{X})^{\#
\mathcal{I}}}([0,\infty))$ be the space of right-continuous functions
with left limits from $[0,\infty)$ to $\mathscr{P}(\mathbb{X})^{\#
\mathcal{I}}$ and $C_{\mathscr{P}(\mathbb{X})^{\#\mathcal
{I}}}([0,\infty))$
the corresponding subset of continuous functions.
The following theorem, which is the main result of the section, shows
that the measure-valued system of Proposition~\ref
{theorem:mv-generator} converges in distribution to a collection of
interacting Fleming--Viot processes. These generalize the celebrated
class of Fleming--Viot diffusions, which take values in the space of
probability measures, to a system of dependent diffusion processes.
See Appendix~\ref{appmA} for a brief review of the essential features.
Here convergence in distribution means weak convergence of the sequence
of distributions induced for each $n$ by $Y^{(n)}(\cdot)$ (as in
Proposition~\ref{theorem:mv-generator}) onto the space $D_{\mathscr
{P}(\mathbb{X})^{\#\mathcal{I}}}([0,\infty))$, to that induced on
the same space by a system of interacting Fleming--Viot diffusions,
with the limiting measure concentrated on $C_{\mathscr{P}(\mathbb
{X})^{\#\mathcal{I}}}([0,\infty))$.

%th4.3 #&#
\begin{theorem}\label{theorem:convergence}
Let $Y^{(n)}(\cdot)=\{Y^{(n)}(t),t\ge0\}$ be as in Proposition \ref
{theorem:mv-generator} with initial distribution $Q_{n}\in(\mathscr
{P}^{n}(\mathbb{X}))^{\#\mathcal{I}}$, and let $Y(\cdot)=\{Y(t),t\ge
0\}$ be a system of interacting Fleming--Viot processes with initial
distribution $Q\in(\mathscr{P}(\mathbb{X}))^{\#\mathcal{I}}$ and
generator defined in Appendix~\textup{\ref{appmA}} by (\ref{F})--(\ref
{gen1}). Assume $\mathbb{X}=[0,1]$, $a(\cdot,\cdot)\equiv m(\cdot
,\cdot)$ and $M^{*}(x,\d y)=\nu_{0}(\d y)$. If additionally $\sigma$
in (\ref{gen1}) is univariate, $\lambda_{n}=\mathrm{O}(n^{2}\#\mathcal{I})$
and $Q_{n}\Rightarrow Q$,
then
\[
Y^{(n)}(\cdot)\Rightarrow Y(\cdot)\qquad\mbox{as }n\rightarrow
\infty
\]
in the sense of convergence in distribution in $C_{\mathscr{P}(\mathbb
{X})^{\#\mathcal{I}}}([0,\infty))$.
\end{theorem}

%%%%%%%%%%%%%%%%%%%%%%%%%%%%%%%%%%%%%%%%%%%%%%%%%%%%%%%%%%%%%%%%%
%%%%%%%%%%%%%%%%%%%%%%%%%%%%%%%%%%%%%%%%%%%%%%%%%%%%%%%%%%%%%%%%%

%s5 #&#
%s5 ###
\section{Algorithms and simulation study}\label{sec:simulation}

In this section we device suitable simulation schemes for the above
constructed systems by means of Markov chain Monte Carlo techniques.
This allows us to explore different economic scenarios and perform
sensitivity analysis on the effects of the model parameters on the
regime changes. Remark~\ref{remark:gibbs sampler} points out that the
discrete representation for a single market can be obtained by means of
Gibbs sampling the joint distribution $q_{n,i}$ in (\ref{joint}). A
similar statement holds for the particle system in a multi market framework.
%
%The Gibbs sampler is special case of the Metropolis-Hastings algorithm
%where the proposal distribution is the conditional distribution of one
%component of the vector whenever we leave unchanged all other
%components. Suppose we have an $n$-dimensional vector, we want to
%sample from the target distribution $p(x_{1},\ldots,x_{n})$, and the
%p(x_{i}|x_{1},\ldots,x_{i-1},x_{i+1},\ldots,x_{n})\nonumber
%are easy to sample from. Then given an initial set of values
%$(x^{0}_{1},\ldots,x^{0}_{n})$, the Gibbs sampler sequentially updates
%all coordinates, namely
%x^{1}_{1}\sim& p(x_{1}|x^{0}_{2},\ldots,x^{0}_{n})\\
%x^{1}_{2}\sim& p(x_{2}|x^{1}_{1},x^{0}_{3},\ldots,x^{0}_{n})\\
%x^{1}_{n}\sim& p(x_{n}|x^{1}_{1},\ldots,x^{1}_{n-1})\\
%x^{2}_{1}\sim& p(x_{1}|x^{1}_{2},\ldots,x^{1}_{n})
%and so on. This algorithm generates an $n$-dimensional Markov chain
%whose stationary distribution is $p(x_{1},\ldots,x_{n})$. Also, each
%component $x_{i}$ has stationary distribution given by the marginal
%$p(x_{i})$. If the updating mechanism is done in a random order, under
%some weak conditions, then reversibility with respect to $p(x_{1},
The particle system in Section~\ref{sec:constr} is such that after a
market $r$ and an item $x_{i}^{r}$ are chosen with probability $\varrho
_{r}$ and $\gamma_{n,i}^{r}$ respectively, a new value for $x_{i}^{r}$
is sampled from
\[
\label{ancora q}
q_{_{D_{n},i}}(\d x^{r}_{i}| \mathcal{I}(-x^{r}_{i}))
\propto\beta_n(x^{r}_i) \alpha_{_{\mathcal{I}(-r)}}(\d
x^{r}_{i})+\sum_{k\ne i}^n\beta_n(x^{r}_k) \delta_{x^{r}_k}(\d x^{r}_{i}),
\]
which selects the next ownership of the vacant share, and all other
items are left unchanged. It is clear that $q_{_{D_{n},i}}$ is the full
conditional distribution of $x_{i}^{r}$ given the current state of the
system. Since the markets, and the particles within the markets, are
updated in random order, it follows immediately that the particle
system is reversible, hence stationary, with respect to (\ref
{joint_particles}).

Algorithm~\ref{alg1} is the random scan Gibbs sampler which generates a sample
path of the particle system with the desired number of markets. Here
we restrict to the case of $\sigma\equiv0$, which implies that the
normalizing constant $\bar q_{_{D_{n},i}}$ is $\theta+n-1$.

\begin{algorithm}
\caption{}\label{alg1}
Initialize; then:
\begin{enumerate}[3.]
\item[1.] select a market $r$ with probability $\varrho_{r}$;
\item[2.] within $r$, select $x^{r}_{i}$ with probability $\gamma^{r}_{n,i}$;
\item[3.] sample $u\sim\operatorname{Unif}(0,1)$;
\item[4.] update $x^{r}_{i}$:
\begin{enumerate}
\item[a.] if $u<\frac{\pi\theta}{\theta+n-1}$, sample
$x^{r}_{i}\sim\nu_{0}$;
\item[b.] if $u>\frac{\theta}{\theta+n-1}$, sample uniformly an
$x^{r}_{k}$, $k\ne i$, within market $r$ and set $x^{r}_{i}=x^{r}_{k}$;
\item[c.] else:
\begin{enumerate}[ii.]
\item[i.] select a market $r'$ with probability $m(r,r')$;
\item[ii.] sample uniformly an $x^{r'}_{j}$ within market $r'$ and
set $x^{r}_{i}=x^{r'}_{j}$;
\end{enumerate}
\end{enumerate}
\item[5.] go back to 1.
\end{enumerate}
\end{algorithm}
%

%re5.1 #&#
\begin{remark}\label{remark:fractionalization}
Note that the fact that updating the whole vector implies sampling from
$n$ different distributions does not lead to an increase in
computational costs if one wants to simulate from the model. Indeed,
acceleration methods such as those illustrated in~\cite{IJ01} can be
easily applied to the present framework.
\end{remark}

As previously mentioned, setting $\sigma\equiv0$, hence $\beta\equiv
1$, as in Algorithm~\ref{alg1}, does not lead to neutrality among firms,
determining instead a competitive advantage of the largest (in terms of
shares) on the smallest. A different choice for $\beta$ allows us to
correct or change arbitrarily this feature. For example, choosing
$\beta(x^{r}_{j^{*}},\mu_{r})=n_{j}^{-1}$, where $n_{j}$ is the
absolute frequency associated with cluster $x^{r}_{j^{*}}$, yields
actual neutrality.
% In general, for $\beta(x^{r}_{k},\mu_{r})$, (\ref{pred2}) and (
% q_{_{D_{n},i}}(\d x^{r}_i|& \mathcal{I}(-x^{r}_{i}))
%& +\theta(1-\pi)\beta_n(x^{r}_i,\mu_{r})\sum_{r'\in\mathcal{I}}m(r,r')
% +\sum_{k\ne i}^n\beta_n(x^{r}_k,\mu_{r}) \delta_{x^{r}_k}(\d x^{r}_i)
%and
%+\sum_{k\ne i}^n\beta_n(x^{r}_k,\mu_{r})\nonumber
Observe also that sampling from (\ref{pred2}), which is composed of
three additive terms, is equivalent to sampling either from
%
%e35 #&#
\begin{equation}\label{alg2-distr1}
\frac{\beta_n(x^{r}_i,\mu_{r})\nu_{0}(\d x^{r}_{i})}
{\int\beta_n(y,\mu_{r})\nu_{0}(\d y)}
\end{equation}
with probability
\[
\frac{\theta\pi}{\bar q_{_{D_{n},i}}}\int\beta_n(y,\mu_{r})\nu
_{0}(\d y),
\]
from
%
%e36 #&#
\begin{equation}\label{alg2-distr2}
\frac{\sum_{r'\in\mathcal{I}}m(r,r')\sum_{j=1}^{n}\beta
_n(x^{r'}_j,\mu_{r})\delta_{x^{r'}_{j}}(\d x^{r}_{i})}
{\sum_{r'\in\mathcal{I}}m(r,r')\sum_{j=1}^{n}\beta_n(x^{r'}_j,\mu_{r})}
\end{equation}
with probability
\[
\frac{\theta(1-\pi)}{\bar q_{_{D_{n},i}}}\sum_{r'\in\mathcal
{I}}m(r,r')\frac{1}{n}\sum_{j=1}^{n}\beta_n(x^{r'}_j,\mu_{r})
\]
or from
%
%e37 #&#
\begin{equation}\label{alg2-distr3}
\frac{\sum_{k\ne i}^n\beta_n(x^{r}_k,\mu_{r}) \delta_{x^{r}_k}(\d x^{r}_i)}
{\sum_{k\ne i}^n\beta_n(x^{r}_k,\mu_{r})}
\end{equation}
with probability
\[
\frac{1}{\bar q_{_{D_{n},i}}}\sum_{k\ne i}^n\beta_n(x^{r}_k,\mu_{r}),
\]
where the normalizing constant $\bar q_{_{D_{n},i}}$ is given by
\begin{eqnarray*}%\label{pred2}
\theta\pi\int\beta_n(x)\nu_{0}(\d x)
+\theta(1-\pi)\sum_{r'\in\mathcal{I}}m(r,r')\int\beta_n(x)\mu
_{r'}(\d x)
+\sum_{k\ne i}^{n}\beta_n(x^{r}_k).
\end{eqnarray*}
Once the functional forms for $\beta$ and $m$ are chosen, computing
$\bar q_{_{D_{n},i}}$ is quite straightforward. If, for example,
$\mathbb{X}=[0,1]$, and the type of an individual admits also
interpretation as index of relative advantage, then one can set $\beta
(x)=x$, and $\bar q_{_{D_{n},i}}$ becomes
\begin{eqnarray*}
\theta\pi{\bar\nu_{0}}
+\theta(1-\pi)\sum_{r'\in\mathcal{I}}m(r,r')\bar{x}^{r'}
+\sum_{k\ne i}^{n}x^{r}_k,
\end{eqnarray*}
where ${\bar\nu_{0}}$ is the mean of $\nu_{0}$, and $\bar{x}^{r'}$
is the average of the components of market $r'$.
To this end, note also that the assumption of $\nu_{0}$ being
non-atomic can be relaxed simplifying the computation.
Algorithm~\ref{alg2} is the extended algorithm for $\beta_n\not\equiv1$.

\begin{algorithm}[t]
\caption{}\label{alg2}
Initialize; then:
\begin{enumerate}[3.]
\item[1.] select a market $r$ with probability $\varrho_{r}$;
\item[2.] within $r$, select $x^{r}_{i}$ with probability $\gamma^{r}_{n,i}$;
\item[3.] sample $u\sim\operatorname{Unif}(0,1)$;
\item[4.] update $x^{r}_{i}$:
\begin{enumerate}[c.]
\item[a.] if $u<\bar q_{_{D_{n},i}}^{-1}\pi\theta\int\beta
_n(y,\mu_{r})\nu_{0}(\d y)$, sample $x^{r}_{i}$ from (\ref{alg2-distr1});
\item[b.] if $u>1-\bar q_{_{D_{n},i}}^{-1}\sum_{k\ne i}^n\beta
_n(x^{r}_k,\mu_{r})$, sample $x^{r}_{i}$ from (\ref{alg2-distr3});
\item[c.] else sample $x^{r}_{i}$ from (\ref{alg2-distr2});
\end{enumerate}
\item[5.] go back to 1.
\end{enumerate}\vspace*{-1pt}
\end{algorithm}
%
%f1 #&#
%f1 ###
\begin{figure}[b]

\includegraphics{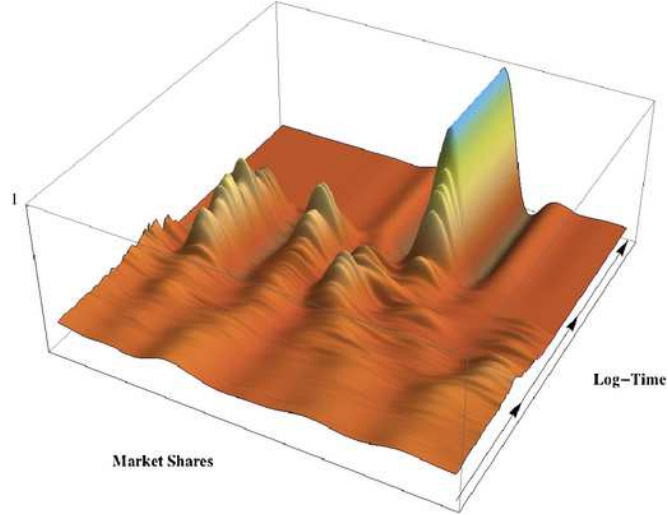}

\caption{High sunk costs progressively transform a perfectly
competitive market into an oligopoly and then into a monopoly.}\label{fig:0a}
\end{figure}

In the following we illustrate how the above algorithms produce
different scenarios where economic regime transitions are caused or
affected by the choice of parameters, which can be structural or
imposed by the policy maker during the observation period. We first
consider a single market and then two interacting markets, and for
simplicity we confine to the use of Algorithm~\ref{alg1}. As a common setting to
all examples we take $\mathbb{X}=[0,1]$, $n=500$, $\nu_{0}$ to be the
probability distribution corresponding to a $\operatorname{Beta}(a,b)$
random variable, with $a,b>0$, with the state space discretized into 15
equally spaced intervals. The number of iterations is $5\times10^{5}$,
of which about 150 are retained at increasing distance. Every figure
below shows the time evolution of the empirical measure of the market,
which describes the concentration of market shares, where time is in
log scale.

Figure~\ref{fig:0a} shows a single market which is in an initial state
of balanced competition among firms, which have similar sizes and
market shares: this can be seen by the flat side closest to the reader.
As time passes, though, the high level of sunk costs, determined by
setting a low $\theta$, is such that exits from the market are not
compensated by the entrance of new firms, and a progressive
concentration occurs. The competitive market first becomes an
oligopoly, shared by no more than three or four competitors, and
eventually a monopoly. Here $\nu_{0}$ corresponds to a $\operatorname
{Beta}(1,1)$ and $\theta=1$.
The fact that the figure shows the market attaining monopoly and
staying there for a time greater than zero could be interpreted as
conflicting with the diffusive nature of the process with positive
(although small) entrance rate of new firms (mutation rate in
population genetics terms). In this respect it is to be kept in mind,
as already mentioned, that the figure is based on observations farther
and farther apart in time. So the picture does not rule out the
possibility of having small temporary deviations from the seeming
fixation at monopoly, which, however, do not alter the long-run overall
qualitative behavior.\looseness=1

%
%f2 #&#
%f2 ###
\begin{figure}

\includegraphics{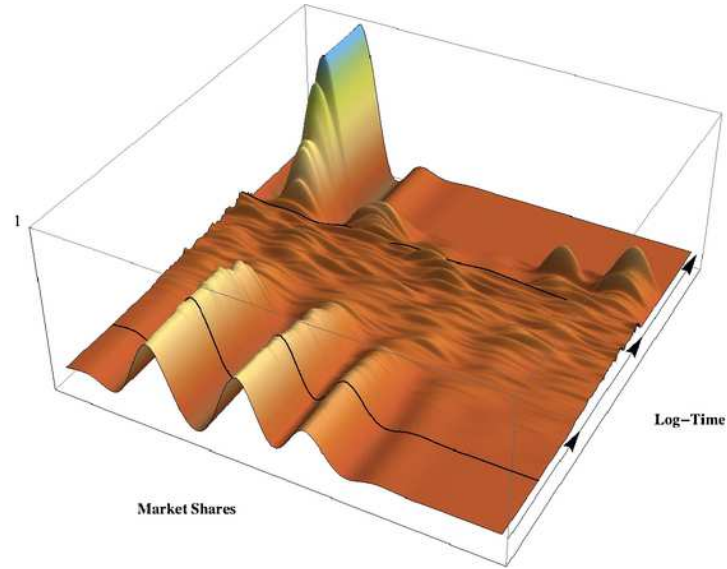}

\caption{An Oligopoly becomes a competitive market after the policy
maker reforms the sector regulation (threshold 1), and concentrates
again after the reform is abolished (threshold~2).}\label{fig:0b}
\end{figure}

In Figure~\ref{fig:0b} we observe a different type of transition. We
initially have an oligopolistic market with three actors. The
structural features of the market are such that the configuration is
initially stable, until the policy maker, in correspondence to the
first black solid line, introduces some new regulation which abates
sunk costs or barriers to entry. Note that in the single market case
the parameter $\theta$ can represent both, since this corresponds to
setting $\pi=1$ in (\ref{alpha}), while in a multiple market
framework we can distinguish the two effects by means of the joint use
of $\theta$ and~$\pi$. Here all parameters are as in Figure~\ref
{fig:0a}, except $\theta$, which is set equal to 1 up to iteration
200, equal to 100 up to iteration $4.5\times10^{4}$ and then equal to 0.
The concentration level progressively decreases and the oligopoly
becomes a competitive market with multiple actors. In correspondence of
the second threshold, namely the second black solid line, there is a
second regulation change in the opposite direction. The market
concentrates again, and, from this point onward, we observe a dynamic
similar to Figure~\ref{fig:0a} (recall that time is in log scale, so
graphics are compressed toward the farthest side). The two thresholds
can represent, for example, the effects of government alternation when
opposite parties have very different political views about a certain sector.

%
%f3 #&#
%f3 ###
\begin{figure}[t]

\includegraphics{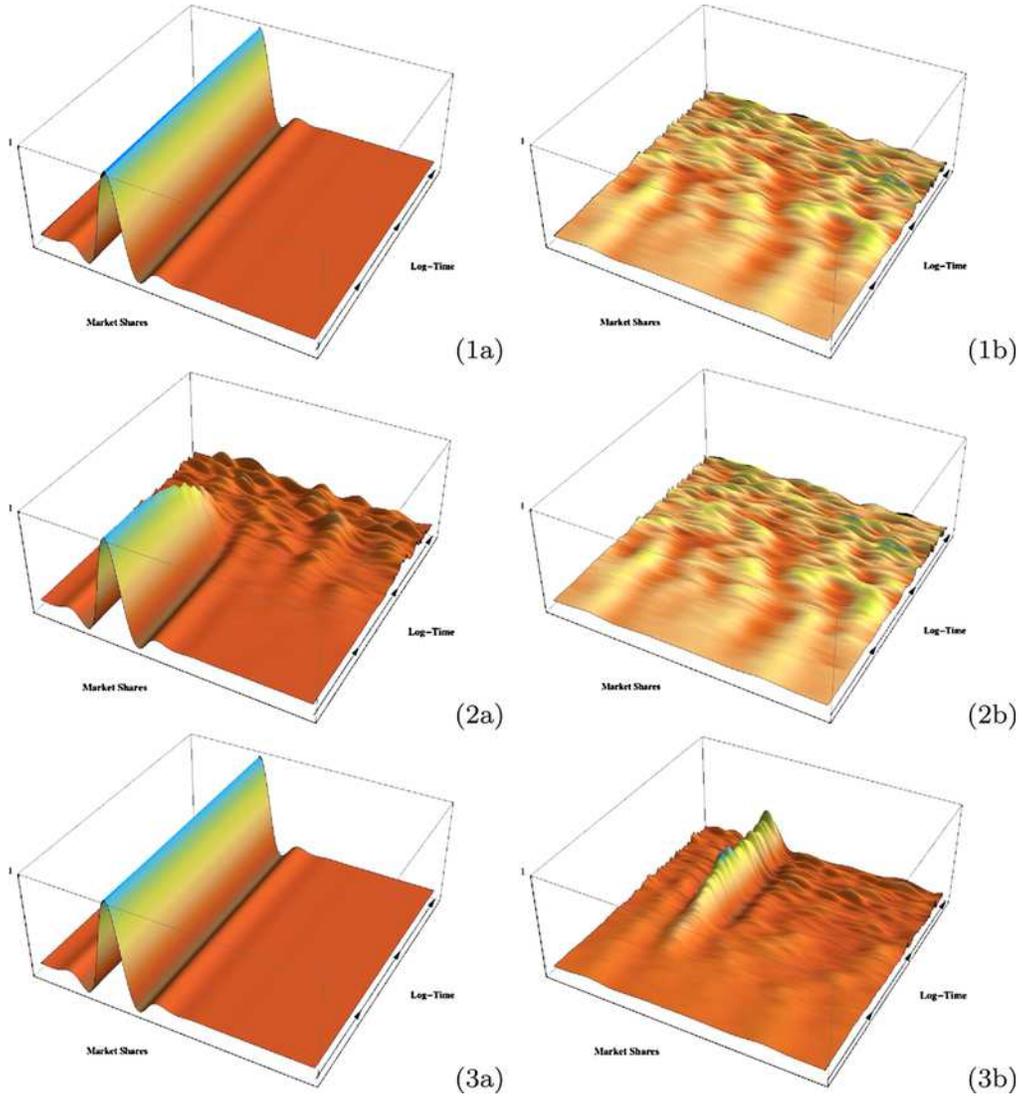}\vspace*{-3pt}

\caption{Effects of parameters' change in interacting monopolistic and
competitive markets. (1a) and (1b) are both closed, hence independent,
markets. (2b) is closed, but (2a) has low barriers to entry ($\pi
\approx0$), and firms from (2b) progressively lower the concentration
in (2a). (3b) has low barriers to entry, so that the monopolist of (3a)
enters the market and conquers a dominant position.}
\label{fig:3}\vspace*{-3pt}
\end{figure}

We now proceed to illustrate some effects of the interaction between
two markets with different structural properties and regulations when
some of
these parameters change. Figure~\ref{fig:3} shows three scenarios
regarding a monopolistic (left) and a competitive market (right).
In all three cases $\nu_{0}$ corresponds to a $\operatorname
{Beta}(1,1)$ for both markets.
Case 1 represents independent markets, due to very high technological
conversion costs or barriers to entry, which is for comparison purposes.
Here $\theta_{a}=0, \theta_{b}=100$ and $\pi_{b}=1$.
In Case 2 the monopolistic market has low barriers to entry, while (2b)
is still closed, and a transition from monopoly to competition occurs.
Here $\theta_{a}=30, \theta_{b}=100$, $\pi_{a}=0.01$, $\pi_{b}=1$.
Case 3 shows the opposite setting, that is, a natural monopoly and a
competitive market with low barriers to entry. The monopolist enters
market (3b) and quickly assumes a dominant position.
Here $\theta_{a}=0, \theta_{b}=100$, $\pi_{b}=0.7$.
Recall, in this respect, the implicit effect due to setting $\beta
\equiv1$, commented upon above.

%
%f4 #&#
%f4 ###
\begin{figure}

\includegraphics{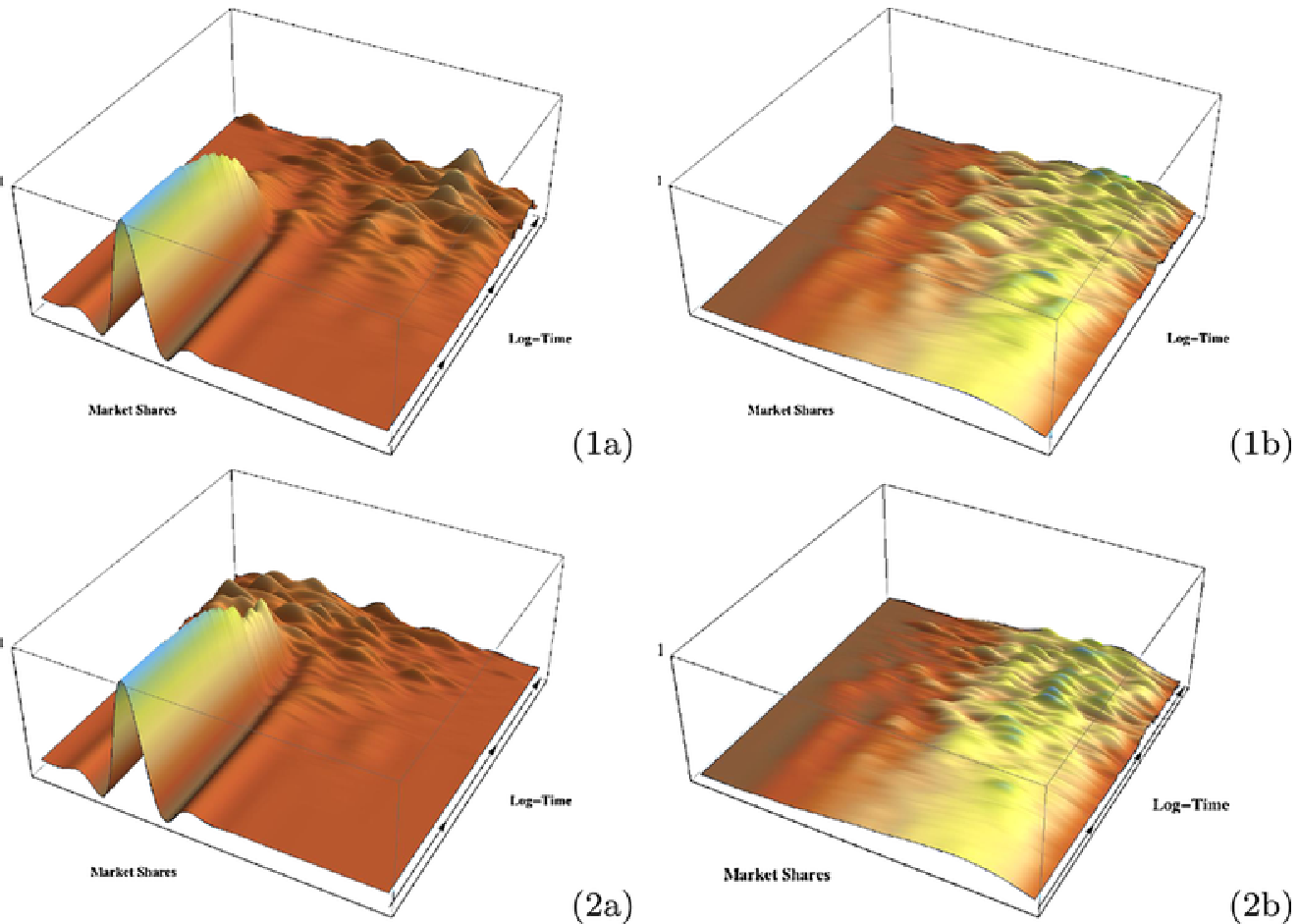}

\caption{Firms in the competitive market (right) are polarized towards
the right half of the state space. (1a)~is a monopoly with high sunk
costs and low barriers to entry, so firms from (1b) enters market (1a).
(2a)~is a monopoly with high barriers to entry and low sunk costs, so
that a transition to a competitive regime occurs independently of
(2b).} \label{fig:9-10}
\end{figure}

Case (2a) in Figure~\ref{fig:3} suggests another point. The
construction of the particle system by means of the hierarchical models
defined in Section~\ref{sec:model} compels us to have the same
centering measure $\nu_{0}$, which generates new firms for all
markets. In particular, this makes it essentially impossible to
establish, by mere inspection of Figure~\ref{fig:3}(2a), whether the
transition is due to new firms or to entrances from (2b). Relaxing this
assumption on $\nu_{0}$ partially invalidates the underlying framework
above, in particular, due to the fact that one loses the symmetry
implied by Lemma~\ref{lemma: P_2n}. Nonetheless the validity of the
particle system is untouched, in that the conditional distributions of
type (\ref{glob-pred2}) are still available, where now $\nu_{0,r}$,
in place of a common $\nu_{0}$, is indexed by $r\in\mathcal{I}$.
This enables us to appreciate the difference between the two above
mentioned effects. If one is willing to give a specific meaning to the
location of the point $x\in\mathbb{X}$ which labels the firm, then
$\nu_{0,r}\ne\nu_{0,r'}$ can model the fact that, say, in two
different sectors, firms are polarized on opposite sides of $\mathbb
{X}$, which, in turn, represents some measurement of a certain
exogenous feature possessed by those firms. Consider a monopoly and a
competitive market, where we now take $\nu_{0,a}$ and $\nu_{0,b}$ to
be the probability measures corresponding to a Beta($2,4$) and a
Beta($4,2$) random variable for the monopolistic and competitive
market, respectively. We are assuming that firms on the left half of
the state space have a certain degree of difference, with respect to
those on the other side, in terms of a certain characteristic. Figure
\ref{fig:9-10} shows the different impact of barriers to entry and
sunk costs on the monopolistic market, due to the joint use of $\pi$
and $\theta$, thus splitting Figure~\ref{fig:3}(2a) into two
different scenarios. The competitive market is composed by firms which
are polarized toward the right half of the state space, meaning, for
example, that they have a high level of a certain feature. Then case 1
of Figure~\ref{fig:9-10} shows the monopoly when sunk costs are high,
but barriers to entry are low, so that the concentration is lowered by
entrance of firms from the other market, rather than from the creation
of new firms from within; while case 2 shows the effects of high
barriers to entry and low sunk costs, so that a transition to a
competitive regime occurs independently of (2b). The parameters for
case 1 are $\theta_{a}=30$, $\theta_{b}=100$, $\pi_{a}=0$, $\pi
_{b}=1$, while for case 2 we have $\theta_{a}=30$, $\theta_{b}=100$,
$\pi_{a}=1$, $\pi_{b}=1$.

%s6 #&#
%s6 ###
\section{Concluding remarks}
In this paper we propose a model for market share dynamics which is
both well founded, from a theoretical point of view, and easy to
implement, from a practical point of view. In illustrating its features
we focus on the impact of changes in market characteristics on the
behaviors of individual firms taking a macroeconomic perspective.
An enrichment of the model could be achieved by incorporating exogenous
information via sets of covariates. This can be done, for example, by
suitably adapting the approach recently undertaken in~\cite{PD10} to
the present framework.
Alternatively, and from an economic viewpoint, more interestingly,
one could modify the model adding a microeconomic understructure:
this would consist of modeling explicitly the individual behavior by
appropriately specifying the function $\beta_n$ at Point (d) in
Section~\ref{sec:constr}, which can account for any desired behavioral
pattern of a single firm depending endogenously on both the status of
all other firms and the market characteristics. This additional layer
would provide a completely explicit micro-foundation of the model,
allowing us to study the effect of richer types of heterogeneous
individual decisions on industry and macroeconomic dynamics through
comparative statics and dynamic sensitivity analysis. These issues of
more economic flavor will be the focus of a forthcoming work.

%apA #&#
\begin{appendix}
%sB #&#
%s7 ###
\section{Background material}\label{appmA}

%The Gibbs sampler is a special case of the Metropolis-Hastings
%algorithm, which in turn belongs to the class of Markov chain Monte
%Carlo procedures. See, e.g.,~\cite{GS90}. These are often applied to
%solve integration and optimization problems in large dimensional
%spaces. Suppose the integral of $f:\mathbb{X}\rightarrow
%$\pi\in\mathscr{P}(\mathbb{X})$ is to be evaluated, and Monte Carlo
%integration turns out to be unfeasible. Markov chain Monte Carlo
%methods provide a way of constructing a stationary Markov chain with $
%first, say, $N$ iterations, and regard the successive
%output from the chain as approximate correlated samples from $\pi$,
%which are then used to approximate $\int f\d\pi$.
%The
%construction of a Gibbs sampler is as follows. Consider a law $\pi=\pi(
%assume that the conditional distributions
% are
%available for every $1\le i\le n$. Then, given an initial set of
%values $(x_{1}^{0},\ldots,x_{n}^{0})$, update iteratively
%x_{1}^{1}&\sim\pi(\d x_{1}|x_{2}^{0},\ldots,x_{n}^{0})\\
%x_{2}^{1}&\sim\pi(\d x_{2}|x_{1}^{1},x_{3}^{0},\ldots,x_{n}^{0})\\
%& \vdots\\
%x_{n}^{1}&\sim\pi(\d x_{n}|x_{1}^{1},\ldots,x_{n-1}^{1})\\
%x_{1}^{2}&\sim\pi(\d x_{1}|x_{2}^{1},\ldots,x_{n}^{1}),
%and so on. Under mild conditions, this routine produces a Markov chain
%with equilibrium law $\pi(\d x_{1},\ldots,\d x_{n})$. The
%above updating rule is known as a \emph{deterministic scan}. If
%instead the components are updated in a random order, called

\subsection*{Basic elements on the Gibbs sampler}

The Gibbs sampler is a special case of the Metropolis--Hastings
algorithm, which, in turn, belongs to the class of Markov chain Monte
Carlo procedures; see, for example,~\cite{GS90}. These are often
applied to solve integration and optimization problems in large
dimensional spaces. Suppose the integral of $f\dvtx \mathbb{X}\rightarrow
\mathbb{R}^{d}$ with respect to
$\pi\in\mathscr{P}(\mathbb{X})$ is to be evaluated, and Monte Carlo
integration turns out to be unfeasible. Markov chain Monte Carlo
methods provide a way of constructing a stationary Markov chain with
$\pi$ as the invariant measure. One can then run the chain, discard
the first, say, $N$ iterations, and regard the successive
output from the chain as approximate correlated samples from $\pi$,
which are then used to approximate $\int f\,\d\pi$.
The
construction of a Gibbs sampler is as follows. Consider a law $\pi=\pi
(\d x_{1},\ldots,\d x_{n})$ defined on $(\mathbb{X}^{n},\mathscr
{X}^{n})$, and
assume that the conditional distributions
\[
\pi(\d x_{i}|x_{1},\ldots,x_{i-1},x_{i+1},\ldots,x_{n})
\]
are
available for every $1\le i\le n$. Then, given an initial set of values
$(x_{1}^{0},\ldots,x_{n}^{0})$, update iteratively
\begin{eqnarray*}%\label{}
x_{1}^{1}&\sim& \pi(\d x_{1}|x_{2}^{0},\ldots,x_{n}^{0}),\\
x_{2}^{1}&\sim& \pi(\d x_{2}|x_{1}^{1},x_{3}^{0},\ldots,x_{n}^{0}),\\
& \vdots&\\
x_{n}^{1}&\sim& \pi(\d x_{n}|x_{1}^{1},\ldots,x_{n-1}^{1}),\\
x_{1}^{2}&\sim& \pi(\d x_{1}|x_{2}^{1},\ldots,x_{n}^{1}),
\end{eqnarray*}
and so on. Under mild conditions, this routine produces a Markov chain
with equilibrium law $\pi(\d x_{1},\ldots,\d x_{n})$. The
above updating rule is known as a \emph{deterministic scan}. If
instead the components are updated in a random order, called \emph
{random scan}, one also gets reversibility with respect to $\pi$.\vspace*{-0.5pt}

\subsection*{Basic elements on Fleming--Viot processes}\vspace*{-0.5pt}

Fleming--Viot processes, introduced in~\cite{FV79}, constitute,
together with Dawson--Watanabe superprocesses, one of the two most
studied classes of probability-measure-valued diffusions, that is,
diffusion processes which take values on the space of probability measures.
A review can be found in~\cite{EK93}.

A Fleming--Viot process can be seen as a generalization of the neutral
diffusion model. This describes the evolution of a vector
$z=(z_{i})_{i\in S}$ representing the relative frequencies of
individual types in an infinite population, where each type is
identified by a point in a space $S$. The process takes values on the simplex\vspace*{-0.5pt}
\[
\Delta_{S}= \biggl\{(z_{i})_{i\in S}\in[0,1]^{S}\dvt z_{i}\ge0, \sum_{i\in
S}z_{i}=1 \biggr\}\\[-0.5pt]
\]
and is characterized by the infinitesimal operator\vspace*{-0.5pt}
\[
L=\frac{1}{2}\sum_{i,j\in S}z_{i}(\delta_{ij}-z_{j})\frac{\partial
^{2}}{\partial z_{i}\,\partial z_{j}}
+\sum_{i\in S}b_{i}(z)\frac{\partial}{\partial z_{i}},\\[-0.5pt]
\]
defined, for example, on the set $C(S)$ of continuous functions on $S$,
if $S$ is compact. Here the first term drives the random genetic drift,
which is the diffusive part of the process, and $b_{i}(z)$ determines
the drift component, with\vspace*{-0.5pt}
\[
b_{i}(z)=\sum_{j\in S,j\ne i}q_{ji}z_{j}-\sum_{j\in S,j\ne i}q_{ij}z_{i}
+z_{i} \biggl(\sum_{j\in S}\sigma_{ij}z_{j}-\sum_{k,l\in S}\sigma
_{kl}z_{k}z_{l} \biggr),\\[-0.5pt]
\]
where $q_{ij}$ is the intensity of a mutation from type $i$ to type $j$
and $\sigma_{ij}=\sigma_{ji}$ is the selection term in a diploid
model. This specification is valid for $S$ finite, which yields the
classical Wright--Fisher diffusion, or countably infinite; see, for
example,~\cite{E81}. Fleming and Viot~\cite{FV79} generalized to the
case of an uncountable type space $S$ by characterizing the
corresponding process, which takes values in the space $\mathscr
{P}(S)$ of Borel probability measures on $S$, endowed with the topology
of weak convergence. Its generator
%can be written
% \mathbb{L}\phi(\mu)=& \frac{1}{2} \int_\mathcal{X}\int_\mathcal{X}\mu(
% {\partial\mu(x)\partial\mu(y) }\\
% &+ \int_\mathcal{X}\mu(\d x) M (\frac{\partial\phi(\mu)}{\partial\mu(
% &+\int_\mathcal{X}\int_\mathcal{X}\mu(\d x)\mu(\d y)(\sigma(x,y)-\la
%where $\partial\phi(\mu)/\partial\mu(x)=\lim_{\varepsilon\rightarrow
%0^+}\varepsilon^{-1}\{\phi(\mu+\varepsilon\delta_x)-\phi(\mu)\}$,
on functions $\phi_{m}(\mu)=F(\la
f_1,\mu\ra,\ldots,\la f_m,\mu\ra)=F(\la\mathbf{f},\mu\ra)$,
where $F\in C^2(\mathbb{R}^m)$,
$f_1,\ldots,f_m$ continuous on $S$ and vanishing at infinity, for
$m\ge1$, and $\la f,\mu\ra=\int f\,\d\mu$,
can be written\vspace*{-0.5pt}
\begin{eqnarray*}
\mathbb{L}\phi(\mu)
&=& \frac{1}{2}\sum_{i,j=1}^{m} (\la f_{i}f_{j},\mu\ra-\la
f_{i},\mu\ra\la f_{j},\mu\ra)F_{z_{i}z_{j}}(\la\mathbf{f},\mu\ra
)\\[-0.5pt]
&&{} +\sum_{i=1}^{m}\la Mf_{i},\mu\ra F_{z_{i}}(\la\mathbf{f},\mu
\ra)+\sum_{i=1}^{m} \bigl(\la(f_{i}\circ\pi)\sigma,\mu^{2}\ra-\la
f_{i},\mu\ra\la\sigma,\mu^{2}\ra\bigr)
F_{z_{i}}(\la\mathbf{f},\mu\ra),\nonumber
\end{eqnarray*}
where $\mu^2$ denotes product measure, $\pi$ is the projection onto
the first coordinate, $M$ is the generator of a Markov process on $S$,
known as the \textit{mutation operator}, $\sigma$ is a non-negative,
bounded, symmetric, Borel measurable functions on $S^{2}$, called
\textit{selection intensity function} and $F_{z_{i}}$ is the
derivative of $F$ with respect to its $i$th argument. Recombination can
also be included in the model.

\subsection*{Interacting Fleming--Viot processes}\label{ssec:FV}
\renewcommand{\theequation}{\arabic{equation}}
\setcounter{equation}{37}

Introduced by~\cite{V90}, and further investigated by~\cite{DGV95}
and~\cite{DG99}, a system of interacting Fleming--Viot processes
extends a Fleming--Viot process to a collection of dependent diffusions
of Fleming--Viot type, whose interaction is modeled as migration of
individuals between subdivided populations.
Following~\cite{DG99}, the model without recombination can be
described as follows. Let the type space be the interval $[0,1]$. Each
component of the system is an element of the set $\mathscr{P}([0,1])$,
denoted $\mu_{r}$ and indexed by a countable set $\mathcal{I}$ of elements
$r,r',\ldots.$ %, so that the state space of the system is $(
%for all $r\in\mathcal{I}$.
%
%Define a second order differential operator on functions on $(
%set of finite measures.
%
%, define
%where
%The second derivative $\frac{\partial^{2}F}{\partial\mu_{r}\partial
%also $\mathcal{A}$ be the algebra of funtions on $(
For $F\dvtx  (\mathscr{P}([0,1]))^{\mathcal{I}}\rightarrow\mathbb{R}$ of
the form
%
%e38 #&#
\begin{equation}\label{F}
F(\mu)=\int_{[0,1]}\cdots\int_{[0,1]}f(x_{1},\ldots,x_{m})\mu
_{r_{1}}(\d x_{1})\cdots\mu_{r_{m}}(\d x_{m})
\end{equation}
with $f\in C([0,1]^{m})$, $(r_{1},\ldots,r_{m})\in(\mathcal
{I})^{m}$, $m\in\mathbb{N}$, the generator of a countable system of
interacting Fleming--Viot processes is
%
%e39 #&#
\begin{eqnarray}\label{gen1}
\mathbb{G} F(\mu)
&=& \sum_{r\in\Omega_N} \biggl\{q\int_{[0,1]} \biggl[\int_{[0,1]}\frac
{\partial F(\mu)}{\partial\mu_{r}}(y)M^{*}(x,\d y)
-\frac{\partial F(\mu)}{\partial\mu_{r}}(x) \biggr]\mu_{r}(\d
x)\nonumber\\
&&\hphantom{\sum_{r\in\Omega_N} \biggl\{}{}+c\sum_{r'\in\Omega_N}a(r,r')\int_{[0,1]}(\mu_{r'}-\mu
_{r})(\d x)\frac{\partial F(\mu)}{\partial
\mu_{r}}(x)\nonumber
\\[-8pt]
\\[-8pt]
&&\hphantom{\sum_{r\in\Omega_N} \biggl\{}{}+d\int_{[0,1]}\int_{[0,1]}\frac{\partial^{2} F(\mu)}{\partial
\mu_{r}\,\partial\mu_{r}}(x,y)Q_{\mu_{r}}(\d x,\d y)\nonumber\\
&&\hphantom{\sum_{r\in\Omega_N} \biggl\{}{}+s\int_{[0,1]}\int_{[0,1]}\int_{[0,1]}\frac{\partial F(\mu
)}{\partial\mu_{r}}(x)\sigma(y,z)\mu_{r}(\d y)Q_{\mu_{r}}(\d x,\d z)
\biggr\},
\nonumber
\end{eqnarray}
where the term $Q_{\mu_{r}}(\d x,\d y)=\mu_{r}(\d x)\delta_{x}(\d
y)-\mu_{r}(\d x)\mu_{r}(\d y)$ drives genetic drift,
$M^{*}(x,\d y)$ is a transition density on $[0,1]\times\mathscr
{B}([0,1])$ modeling mutation, $\mathscr{B}([0,1])$ is the Borel
sigma algebra on $[0,1]$, $a(\cdot,\cdot)$ on $\mathcal{I}\times
\mathcal{I}$ such that $ a(r,r')\in[0,1]$ and $\sum_{r}^{}a(r,r')=1$
is a transition kernel modeling migration and $\sigma(\cdot,\cdot)$
is a bounded symmetric selection intensity function on $[0,1]^{2}$. The
non-negative reals $q,c,d,s$ represent, respectively, the rate of
mutation, immigration, resampling and selection.
Let the mutation operator be
%
%e40 #&#
\begin{equation}\label{mut-op0}
Mf(z)=\int[f(y)-f(z) ]M^{*}(x,\d y),
\qquad f\in B(\mathbb{X})
\end{equation}
and the migration operator be
%
%e41 #&#
\begin{equation}\label{mig-op0}
G^{r'}f(z)=\int[f(y)-f(z) ]\mu_{r'}(\d y),
\qquad f\in B(\mathbb{X}),
\end{equation}
for $r'\in\mathcal{I}$.
Using this notation,
% \mathbb{G} F(\mu)
% =&\sum_{r\in\Omega_N} \{q\int_{[0,1]}\mu_{r}(\d x)M (\frac{\partial F(
% &+c\sum_{r'\in\Omega_N}a(r,r')\int_{[0,1]}\mu_{r}(\d x)G^{r'} (\frac{
%% &+r\int_{[0,1]}\int_{[0,1]}\int_{[0,1]}\alpha(y_{2},y_{3})\mu_{r}(\d
%y_{1})\mu_{r}(\d y_{2})\mu_{r}(\d y_{3})\tilde R (\frac{\partial F(
% &+d\int_{[0,1]}\int_{[0,1]}\mu_{r}(\d x) (\delta_{x}(\d y)-\mu_{r}(\d
%y) )\frac{\partial^{2} F(\mu)}{\partial\mu_{r}\partial\mu_{r}}(x,y)
% &+s\int_{[0,1]}\int_{[0,1]}\mu_{r}(\d x)\mu_{r}(\d y) (\sigma(x,y)-
% \end{eqnarray*}
%where $\la f,\mu\ra=\int f\d\mu$. Furthermore,
and when $F$ is as in (\ref{F}), (\ref{gen1}) can be written
%
%e42 #&#
\begin{eqnarray}\label{gen2}
\mathbb{G} F(\mu)&=&\sum_{r\in\Omega_N} \Biggl\{ q\sum_{i=1}^{m}\int
_{[0,1]}\cdots\int_{[0,1]} M_{r_{i}}f\,\d\mu_{r_{1}}\cdots\,\d\mu
_{r_{m}}\nonumber\\
&&\hphantom{\sum_{r\in\Omega_N} \Biggl\{}{}+c\sum_{r'\in\Omega_N}a(r,r')\sum_{i=1}^{m}\int_{[0,1]}\cdots
\int_{[0,1]} G_{r_{i}}^{r'}f\,\d\mu_{r_{1}}\cdots\,\d\mu_{r_{m}}\nonumber\\
&&\hphantom{\sum_{r\in\Omega_N} \Biggl\{}{}+d\sum_{i=1}^{m}\sum_{k\ne i}^{m}\int_{[0,1]}\cdots\int
_{[0,1]} ( \Phi_{r_{k,i}}f-f )\,\d\mu_{r_{1}}\cdots\,\d\mu
_{r_{m}}\\
&&\hphantom{\sum_{r\in\Omega_N} \Biggl\{}{}+s\sum_{i=1}^{m}\int_{[0,1]}\cdots\int_{[0,1]} (\sigma
_{r_{i},m+1}(\cdot,\cdot)f \nonumber\\
&&\hphantom{\sum_{r\in\Omega_N} \Biggl\{{}+s\sum_{i=1}^{m}\int_{[0,1]}\cdots\int_{[0,1]} (}{} - \sigma_{m+1,m+2}(\cdot,\cdot)f )\,\d\mu_{r_{1}}\cdots\,\d\mu
_{r_{m}}\,\d\mu_{r}\,\d\mu_{r}
\Biggr\},
\nonumber
\end{eqnarray}
where $M_{j}$ and $G_{j}^{r'}$ are $M$ and $G^{r'}$
applied to the $j$th coordinate of $f$, $r_{i}$ is as in Proposition~\ref{theorem:particle-gen},
$r_{k,i}$ as in (\ref{xi-oper2}) and
$\Phi_{hj}$ as in (\ref{Phi}). When $\mathcal{I}$ is single-valued,
(\ref{gen2}) simplifies to
\begin{eqnarray*}%\label{DKgen}
\mathbb{G} F(\mu)&=& q\sum_{i=1}^{m}\la M_{i}f,\mu^{m}\ra
+d\sum_{i=1}^{m}\sum_{k\ne i}^{m}\la\Phi_{ki}f-f,\mu^{m}\ra\\
&&{}+s\sum_{i=1}^{m} \bigl(\la\sigma_{i,m+1}(\cdot,\cdot)f,\mu
^{m+1}\ra-\la\sigma_{m+1,m+2}(\cdot,\cdot)f,\mu^{m+2}\ra
\bigr),\nonumber
\end{eqnarray*}
which is the generator of a Fleming--Viot process with selection with
$F(\mu)=\la f,\mu^{m}\ra$, $f\in C([0,1])$.

%sC #&#
%s8 ###
\section{Proofs}\label{appmB}
\setcounter{equation}{42}
%From (\ref{MDP}) together with Proposition 1 and Corollary 3.2 and
%3.2' in~\cite{A74}, we have that
%Y_{1}| x_{1},\ldots,x_{n}\sim& \frac{\alpha_{x_{1},\ldots,x_{n}}(\cdot
%)}{\alpha_{x_{1},\ldots,x_{n}}(\mathbb{X})}
%where $\alpha_{x_{1},\ldots,x_{n}}(\mathbb{X})$ is the total mass of $
%Y_{n}| y_{1},\ldots,y_{n-1},x_{1},\ldots,x_{n}\sim&
%Hence, in view of (\ref{P_n}), the joint law of $(Y_{1},\ldots,Y_{n})$
%conditional on $(X_{1},\ldots,X_{n})$ is
%%p_{n}(\d y_{1},\ldots,\d y_{n}| x_{1},\ldots,x_{n})
%Denoting (\ref{eq: Polya-urn}) by $\mathcal{M}^{\alpha}_{n}(\d
%x_{n}|x_{1},\ldots,x_{n-1})$, it can be easily seen that (\ref{eq:
%alpha_x}) implies
%=\frac{\alpha(\d y)+\delta_{x}(\d y)}{\alpha(\mathbb{X})+1}
%=\mathcal{M}_{2}^{\alpha}(\d y|x)
%and, more generally,
%=\mathcal{M}_{2n}^{\alpha}(\d y_{1},\ldots,\d y_{n}|x_{1},\ldots,x_{n})
%from which
%=\mathcal{M}_{2n}^{\alpha}(\d x_{1},\ldots,\d x_{n},\d y_{1},\ldots,\d
%y_{n}).
%Since the right hand side is symmetric in its arguments (cf. (
\begin{pf*}{Proof of Proposition~\ref{theorem:particle-gen}}
%%\label{}
The infinitesimal generator of the $\mathbb{X}^{n}$-valued process
described at the beginning of Section~\ref{sec:constr} can be written,
for any $f\in B(\mathbb{X}^{n})$, as
%
%e43 #&#
\begin{equation}\label{gen1col}
A_{n}f(\mathbf{x})=\lambda_{n}\sum_{i=1}^n\gamma_{n,i}\int[f(\eta
_i(\mathbf{x}|y))-f(\mathbf{x}) ]q_{n,i}\bigl(\d y|\mathbf{x}_{(-i)}\bigr),
\end{equation}
where $q_{n,i}(\d y|\mathbf{x}_{(-i)})$ is (\ref{pred}) and $\eta
_{i}$ is
as in (\ref{eta}).
%defined by
Within the multi-market framework, (\ref
{gen1col}) is the generator of the process for the configuration of
market $r$, say, conditionally on all markets $r'\in\mathcal{I}$,
$r'\ne r$, and can be written
%
%e44 #&#
\begin{equation}\label{gen1col1b}
A_{D_{n}}f(\mathbf{x}^{r}| \mathcal{I}(-r))=\lambda_{n}\sum
_{i=1}^n\gamma^{r}_{n,i}\int[f(\eta_i(\mathbf{x}^{r}|y))-f(\mathbf
{x}^{r}) ] q_{_{D_{n},i}}(\d y|\mathcal{I}(-x^{r}_{i})),
\end{equation}
where $\mathcal{I}(-r)$ and $\mathcal{I}(-x^{r}_{i})$ are as in (\ref
{eq: Xi-xi}) and (\ref{eq: Xi-x_i}), $\gamma_{n,i}^{r}$ are the
market-specific removal probabilities and $ q_{_{D_{n},i}}(\d
y|\mathcal{I}(-x^{r}_{i}))$ is (\ref{glob-pred2}).
%In order to write the generator of the entire system, some additional
%notation is needed. For a given operator $H$ on functions $g\in B(
%n$, let
%For $h\in\mathbb{N}$ and $\mathbf{x}\in[0,1]^{h}$ such that $x_{k}\in(
%operator on functions $g\in B([0,1]^{h})$ which maps $g(\mathbf{x})
%B([0,1]^{h})$. Denote with
Then the generator for the whole particle system, for every $f\in
B(\mathbb{X}^{D_{n}})$, is
%
%e45 #&#
\begin{equation}\label{gen1col2}
A_{D_{n}}f(\mathbf{x})=\lambda_{n}\sum_{r\in\mathcal{I}}\varrho
_{r}\sum_{i=1}^n\gamma^{r}_{n,i}\int[f(\eta_{r_{i}}(\mathbf
{x}|y))-f(\mathbf{x}) ] q_{_{D_{n},i}}(\d y|\mathcal{I}(-x^{r}_{i})),
\end{equation}
where $\eta_{r_{i}}$ is as in (\ref{xi-oper}).
%where $\eta$ is (\ref{eta}), $\eta_{r_{i}}(\mathbf{x}|y)$ is (
%q_{_{D_{n},i}}(\d y|\mathcal{I}(-x^{r}_{i}))$ is (\ref{glob-pred2}).
%The operator (\ref{gen1col2}) is the generator of a $
%At rate $\lambda_{n}$ a transition occurs. At each transition, a
%market $r$ is randomly selected with probability $\varrho_{r}$ and a
%component of $(x^{r}_{1},\ldots,x^{r}_{n})$ is randomly selected with
%probability $\gamma^{r}_{n,i}$. The next state is obtained by setting
%all components of the system different from $x^{r}_{i}$ equal to their
%previous state, and by sampling a new value for $x^{r}_{i}$ from (
Setting now $\beta_n$ as in (\ref{sel_funct}), (\ref{gen1col2}) becomes
%
%e46 #&#
\begin{eqnarray}\label{gen3}
A_{D_{n}}f(\mathbf{x})
&=&\sum_{r\in\mathcal{I}} \Biggl\{\sum_{i=1}^nC_{n,r,i}\int[f(\eta_{r_{i}}
(\mathbf{x}|y))-f(\mathbf{x}) ] \biggl(1+\frac{2\sigma(y)}{n} \biggr)\alpha
_{_{\mathcal{I}(-r)}}(\d y)\nonumber\\
&&\hphantom{\sum_{r\in\mathcal{I}} \Biggl\{}{}+\sum_{1\le k\ne i\le n}C_{n,r,i} [f(\eta_{r_{i}}(\mathbf
{x}|x^{r}_k))-f(\mathbf{x}) ]\\
&&\hphantom{\sum_{r\in\mathcal{I}} \Biggl\{}{}+\frac{1}{n}\sum_{1\le k\ne i\le n}C_{n,r,i} \sigma(x^{r}_k)
[f(\eta_{r_{i}}(\mathbf{x}|x^{r}_k))-f(\mathbf{x}) ] \Biggr\},\nonumber
\end{eqnarray}
with $C_{n,r,i}$ as in (\ref{constant}).
%Specifying $\alpha_{_{\mathcal{I}(-r)}}$
%Let $\theta>0$, $\pi\in[0,1]$, and let $\nu_{0}$ be a non atomic
%probability measure on $(\mathbb{X},\mathscr{X})$. Then set
%where $\nu_{1}$ is the probability measure
%=\sum_{r'\in\mathcal{I}}m(r,r')\mu_{r'}(\d y).
%Here $\mu_{r'}\in\mathscr{P}(\mathbb{X})$ is the empirical measure of
%the share distribution in market $r'$, namely
%and $m(r,r'):\mathcal{I}\times\mathcal{I}\rightarrow[0,1]$ is such that
%m(r,r)=0,
%
%Hence $\alpha_{_{\mathcal{I}(-r)}}$ is
%as in (\ref{alpha}), (\ref{xi-empirical}) and (\ref{m}),
%with total mass $\alpha_{_{\mathcal{I}(-r)}}(\mathbb{X})=\theta$, by
%means of an analogous argument to that used in Lemma~\ref{lemma: P_2n}
%it can be shown that
%(\ref{alpha}) satisfies (\ref{eq: marginal-requirement2}).
Substituting (\ref{alpha}) in (\ref{gen3}) yields %, for $f\in
%B(\mathbb{X}^{D_{n}})$,
%
%e47 #&#
\begin{eqnarray}\label{genallcol}
A_{D_{n}}f(\mathbf{x})
& =& \sum_{r\in\mathcal{I}} \Biggl\{\theta\pi\sum_{i=1}^{n}C_{n,r,i}
\int
[f(\eta_{r_{i}}(\mathbf{x}|y)) - f(\mathbf{x}) ]\biggl(1 + \frac{\sigma
(y)}{n}\biggr)\nu_{0}(\d y)\nonumber\\
&&\hphantom{\sum_{r\in\mathcal{I}} \Biggl\{}{}+\theta(1-\pi)\sum_{r'}m(r,r')\sum_{1\le j\ne i\leq
n}C_{n,r,i}\int[f(\eta_{r_{i}}(\mathbf{x}|y))-f(\mathbf{x})
]\nonumber\\
&&\hphantom{\sum_{r\in\mathcal{I}} \Biggl\{ {}+\;}{} \times\biggl(1+\frac{\sigma(y)}{n}\biggr)\mu_{r'}(\d y)\\
&&\hphantom{\sum_{r\in\mathcal{I}} \Biggl\{}{}+\sum_{1\le k\ne i\leq n}C_{n,r,i}
[f(\eta_{r_{i}}(\mathbf
{x}|x^{r}_k))-f(\mathbf{x}) ]\nonumber\\
&&\hphantom{\sum_{r\in\mathcal{I}} \Biggl\{}{}+\frac{1}{n}\sum_{1\le k\ne i\le n}C_{n,r,i} \sigma(x^{r}_k)
[f(\eta_{r_{i}}(\mathbf{x}|x^{r}_k))-f(\mathbf{x}) ] \Biggr\}.\nonumber
\end{eqnarray}
By means of (\ref{mut-op}) and (\ref{mig-op}), with $M_{i}^{n}f$ and
$G_{i}^{n,r'}f$ denoting, respectively, $M^{n}$ and $G^{n,r'}$, applied
to the $i$th coordinate of $f$, and $M_{r_{i}}^{n}f$ and
$G_{r_{i}}^{n,r'}f$ interpreted according to (\ref{xi-oper}), (\ref
{genallcol}) can be written
\begin{eqnarray*}%\label{genallcol2}
A_{D_{n}}f(\mathbf{x})
&=&\sum_{r\in\mathcal{I}} \Biggl\{\theta\pi\sum_{i=1}^{n}C_{n,r,i}
M_{r_{i}}^{n}f(\mathbf{x})\\
&&\hphantom{\sum_{r\in\mathcal{I}} \Biggl\{}{}+\theta(1-\pi)\sum_{r'}m(r,r')\sum_{i=1}^{n}C_{n,r,i}
G_{r_{i}}^{n,r'}f(\mathbf{x})\nonumber\\
&&\hphantom{\sum_{r\in\mathcal{I}} \Biggl\{}{}+\sum_{1\le k\ne i\le n}C_{n,r,i}
[f(\eta_{r_{i}}(\mathbf
{x}|x^{r}_k))-f(\mathbf{x}) ]\nonumber\\
&&\hphantom{\sum_{r\in\mathcal{I}} \Biggl\{}{}+\frac{1}{n}\sum_{1\le k\ne i\le n}C_{n,r,i} \sigma(x^{r}_k)
[f(\eta_{r_{i}}(\mathbf{x}|x^{r}_k))-f(\mathbf{x}) ] \Biggr\}.\nonumber
\end{eqnarray*}
\upqed
\end{pf*}
\begin{pf*}{Proof of Proposition~\ref{theorem:mv-generator}}
%%\label{}
For $k\le n$, let $n_{[k]}$ be as in (\ref{desc-factorial}), and
define the probability measure
%
%e48 #&#
\begin{equation}\label{mu-(D_k)}
\mu^{(D_{k})}=\prod_{r\in\mathcal{I}}\frac{1}{n_{[k]}}\sum_{1\le
i_{r,1}\ne\cdots\ne i_{r,k}\le n}\delta_{(x^{r}_{i_{r,1}},\ldots
,x^{r}_{i_{r,k}})},
\end{equation}
where $D_{k}$ is as in (\ref{eq: system-dimension}).
Define also
\[
\phi_{D_{k}}(\mu)=\bigl\la f,\mu^{(D_{k})}\bigr\ra,
\qquad
f\in B(\mathbb{X}^{D_{k}})
\]
and
%
%e49 #&#
\begin{equation}\label{Aphi}
\mathbb{A}_{D_{n}}\phi_{D_{k}}(\mu)=\bigl\la A_{D_{n}}f,\mu^{(D_{k})}\bigr\ra,
\end{equation}
where $\la f,\mu\ra=\int f\,\d\mu$.
%In a single market framework, the above objects would simplify to
% \mu^{(k)}=\frac{1}{n_{[k]}}\sum_{1\le i_{1}\ne\cdots\ne i_{k}\le n}
%$\la f,\mu^{(k)}\ra$ and $\mathbb{A}_{n}\phi_{k}(\mu)=\la A_{n}f,
Then $\mathbb{A}_{D_{n}}\phi_{D_{n}}(\mu)$ is the generator of the
$(\mathscr{P}(\mathbb{X}))^{\#\mathcal{I}}$-valued system (\ref
{Yn}), which from (\ref{genallcol2}), letting $f\in B(\mathbb
{X}^{D_{n}})$ in (\ref{Aphi}), can be written
%
%e50 #&#
\begin{eqnarray}\label{mv-gen}
\mathbb{A}_{D_{n}}\phi_{D_{n}}(\mu)
&=&\sum_{r\in\mathcal{I}} \Biggl[\theta\pi\sum_{i=1}^nC_{n,r,i}\la
M_{r_{i}}^{n}f,\mu^{(D_{n})}\ra\nonumber\\
&&\hphantom{\sum_{r\in\mathcal{I}} \Biggl[}{}+\theta(1-\pi)\sum_{r'}m(r,r')\sum_{i=1}^nC_{n,r,i} \bigl\la
G_{r_{i}}^{n,r'}f,\mu^{(D_{n})}\bigr\ra\nonumber
\\[-8pt]
\\[-8pt]
&&\hphantom{\sum_{r\in\mathcal{I}} \Biggl[}{}+\sum_{1\le k\ne i\le n}C_{n,r,i}\bigl\la\Phi_{r_{k,i}}f-f,\mu
^{(D_{n})}\bigr\ra\nonumber\\
&&\hphantom{\sum_{r\in\mathcal{I}} \Biggl[}{}+\frac{1}{n}\sum_{1\le k\ne i\le n}C_{n,r,i}\bigl\la\sigma
_{r_{k}}(\cdot)(\Phi_{r_{k,i}}f-f),\mu^{(D_{n})}\bigr\ra\Biggr],
\nonumber
\end{eqnarray}
where $\sigma_{r_{k}}(\cdot)$ denotes $\sigma(x^{r}_{k})$ and $\Phi
_{ki}$ is as in (\ref{Phi}).
% and
% note that for $f\in B(\mathbb{X}^{n})$, $k>n$ implies $\Phi_{ki}f=f$,
%and interpret $\Phi_{r_{k,i}}$ according to (\ref{xi-oper2}).
Note now that for $f\in B(\mathbb{X}^{m})$, $m\le D_{n}$, we have
\begin{eqnarray*}%\label{}
M_{r_{i}}^{n}f=f,
\qquad
G_{r_{i}}^{n,r'}f=f,
\qquad
\Phi_{r_{k,i}}f=f,
\qquad\mbox{if } i>m
\end{eqnarray*}
and
\[
\bigl\la\Phi_{r_{k,i}}f,\mu^{(m)}\bigr\ra=\bigl\la f,\mu^{(m)}\bigr\ra,
\qquad
i\le m, m+1\le k\le n.
\]
Given (\ref{mu-(k-xi)}) and (\ref{mu-(m)}), it follows that when
$f\in B(\mathbb{X}^{m})$, $m\le D_{n}$, (\ref{mv-gen}) can be written
\begin{eqnarray*}
\mathbb{A}_{D_{n}}\phi_{m}(\mu)
&=&\sum_{r\in\mathcal{I}} \Biggl\{\theta\pi\sum_{i=1}^mC_{n,r,i}\bigl\la
M_{r_{i}}^{n}f,\mu^{(m)}\bigr\ra\\
%&+\theta\pi(n-k_{r})\sum_{i=1}^mC_{n,r,i}\la M_{r_{i}}^{n,m+1}f,
&&\hphantom{\sum_{r\in\mathcal{I}} \Biggl\{}{}+\theta(1-\pi)\sum_{r'}a(r,r')\sum_{i=1}^mC_{n,r,i} \bigl\la
G_{r_{i}}^{n,r'}f,\mu^{(m)}\bigr\ra\nonumber\\
% &+\frac{\theta(1-\pi)(n-k_{r})}{n-1}\sum_{r'}a(r,r')\sum_{i=1}^m
%C_{n,r,i}\la G_{r_{i}}^{n,m+1,r'}f,\mu^{(m+1)}\ra\nonumber\\
&&\hphantom{\sum_{r\in\mathcal{I}} \Biggl\{}{}+\sum_{1\le k\ne i\le m}C_{n,r,i}\bigl\la\Phi_{r_{k,i}}f- f,\mu
^{(m)}\bigr\ra\nonumber\\
&&\hphantom{\sum_{r\in\mathcal{I}} \Biggl\{}{}+\frac{1}{n}\sum_{i=1}^m\sum_{k\ne i}^{k_{r}}C_{n,r,i}\bigl\la\sigma
_{r_{k}}(\cdot)(\Phi_{r_{k,i}}f-f),\mu^{(m)}\bigr\ra\nonumber\\
&&\hphantom{\sum_{r\in\mathcal{I}} \Biggl\{}{}+\frac{n-k_{r}}{n}\sum_{i=1}^mC_{n,r,i}\bigl\la\sigma_{m+1}(\cdot
)(\Phi_{m+1,r_{i}}f-f),\mu^{(m)}\mu_{r}\bigr\ra\nonumber
%&+\frac{n-k_{r}}{n}\sum_{i=1}^m\sum_{k\ne i}^{k_{r}}C_{n,r,i}\la
%&+\frac{2(n-k_{r})(n-k_{r}-1)}{n(n-1)}\sum_{i=1}^mC_{n,r,i}\la
\Biggr\}.\nonumber
\end{eqnarray*}
\upqed
\end{pf*}
\begin{pf*}{Proof of Theorem \protect\ref{theorem:convergence}}%
For $f\in B(\mathbb{X}^{k})$, $k\ge1$, let $\|f\|=\sup_{x\in\mathbb
{X}^{k}}|f(x)|$.
Observe that (\ref{mut-op}) and (\ref{mig-op}) converge uniformly,
respectively to (\ref{mut-op0}) and (\ref{mig-op0}), as $n$ tends to
infinity, implying
\begin{eqnarray*}
\bigl\|\bigl\la M^{n}_{r_{i}}f,\mu^{(m)}\bigr\ra-\bigl\la M_{r_{i}}f,\mu^{(m)}\bigr\ra\bigr\|&\rightarrow0&,
\qquad f\in B(\mathbb{X}^{m}),
\\
\bigl\|\bigl\la G^{n,r'}_{r_{i}}f,\mu^{(m)}\bigr\ra-\bigl\la G^{r'}_{r_{i}}f,\mu
^{(m)}\bigr\ra\bigr\|&\rightarrow0&,
\qquad f\in B(\mathbb{X}^{m}).
\end{eqnarray*}
Here the supremum norm is intended with respect to the vector $x\in\mathbb
{X}^{m}$ of atoms in $\mu^{(m)}$, with $\mu^{(m)}$ as in (\ref{mu-(m)}).
Let now $\mu^{(k_{r})}_{r}$ be as in (\ref{mu-(k-xi)}), so that $\mu
_{r}=n^{-1}\sum_{i=1}^{n}\delta_{x^{r}_{i}}$.
Then it is easy to check that
\[
\bigl\|\bigl\la f,\mu^{(k_{r})}_{r}\bigr\ra-\bigl\la f,\mu_{r}^{k_{r}}\bigr\ra\bigr\|\rightarrow0,
\qquad f\in B(\mathbb{X}^{k_{r}}),\nonumber
\]
as $n\rightarrow\infty$, where $\mu^{k_{r}}$ denotes a $k_{r}$-fold
product measure $\mu_{r}\times\cdots\times\mu_{r}$, and that
%Letting also $\mu^{(m)}$ be as in (\ref{mu-(m)}), we have that
%
\[
\bigl\|\bigl\la f,\mu^{(m)}\bigr\ra-\la f,\mu^{\times m}\ra\bigr\|\rightarrow0,
\qquad f\in B(\mathbb{X}^{m}),
\]
as $n\rightarrow\infty$, where we have denoted
\[
\mu^{\times m}=\prod_{r\in\mathcal{I}}\mu_{r}^{k_{r}}.
\]
We also have, from (\ref{constant}) $C_{n,r,i}=\lambda_{n} \varrho
_{r} \gamma_{n,i}^{r}/\bar q_{_{D_{n},i}}$,
where $\lambda_{n}$ is the Poisson rate driving the holding times,
$\varrho_{r}=\mathrm{O}(\#\mathcal{I}^{-1})$ and $\gamma^{r}_{n,i}=\mathrm{O}(n^{-1})$
are the probability of choosing market $r$ and $x^{r}_{i}$ respectively
during the update, and $\bar q_{_{D_{n},i}}=\mathrm{O}(n)$ is the normalizing
constant of (\ref{pred2}). Then choosing $\lambda_{n}=\mathrm{O}(n
D_{n})=\mathrm{O}(n^{2}\#\mathcal{I})$
implies $C_{n,r,i}\rightarrow1$ as $n\rightarrow\infty$.
Finally, let $\varphi_{m}\in B(\mathscr{P}(\mathbb{X}^{m}))$ be
%
%e52 #&#
\begin{eqnarray}\label{varphi}
\varphi_{m}(\mu)=\la f,\mu^{\times m}\ra
=\int_{[0,1]}\cdots\int_{[0,1]}f(x_{1},\ldots,x_{m})\mu_{r_{1}}(\d
x_{1})\cdots\mu_{r_{m}}(\d x_{m})
\end{eqnarray}
for any sequence $(r_{1},\ldots,r_{m})\in\mathcal{I}^{m}$.
Then it can be checked that (\ref{mv-gen-(m)}) converges, as $n$ tends
to infinity, to
%
%e53 #&#
\begin{eqnarray*}%\label{mv-gen-m}
\mathbb{A}\varphi_{m}(\mu) &=&\sum_{r\in\mathcal{I}} \Biggl[\theta\pi
\sum_{i=1}^m\la M_{r_{i}}f,\mu^{\times m}\ra+\theta(1-\pi)
\sum_{r'\in\mathcal{I}}m(r,r')\sum_{i=1}^m \la G_{r_{i}}^{r'}f,\mu
^{\times m}\ra\\
&&\hphantom{\sum_{r\in\mathcal{I}} \Biggl[}{} +\sum_{1\le k\ne i\le m}\la\Phi_{r_{k,i}}f-f,\mu^{\times m}\ra
+\sum_{i=1}^m\la
\sigma_{r_{i}}(\cdot)f- \sigma_{m+1}(\cdot) f,\mu^{\times m}\mu
_{r}\ra
\Biggr]
\end{eqnarray*}
which, in turn, implies
\[
\|\mathbb{A}_{D_{n}}\phi_{m}(\mu)-\mathbb{A}\varphi_{m}(\mu) \|\longrightarrow0
\qquad\mbox{as } n\rightarrow\infty.
\]
Using (\ref{varphi}), and letting $\mathbb{X}=[0,1]$, $\mathbb
{A}\varphi_{m}(\mu)$ can be written
%
%e54 #&#
\begin{eqnarray}\label{mv-gen-m3}
\mathbb{A}\varphi_{m}(\mu) &=&\sum_{r\in\mathcal{I}} \Biggl[\theta\pi
\sum_{i=1}^m\int_{[0,1]}
\cdots\int_{[0,1]} M_{r_{i}}f \,\d\mu_{r_{1}}\cdots\,\d\mu
_{r_{m}}\nonumber\\
&&  \hphantom{\sum_{r\in\mathcal{I}} \Biggl[}{}+\theta(1-\pi)\sum_{r'\in\mathcal{I}}m(r,r')\sum_{i=1}^m
\int_{[0,1]}\cdots\int_{[0,1]} G_{r_{i}}^{r'}f\,\d\mu_{r_{1}}\cdots\,\d\mu_{r_{m}}\nonumber
\\[-8pt]
\\[-8pt]
&&  \hphantom{\sum_{r\in\mathcal{I}} \Biggl[}{}+\sum_{1\le k\ne i\le m}\int_{[0,1]}\cdots\int_{[0,1]}
(\Phi_{r_{k,i}}f-f )\,\d\mu_{r_{1}}\cdots\,\d\mu_{r_{m}}\nonumber\\
&&  \hphantom{\sum_{r\in\mathcal{I}} \Biggl[}{}+\sum_{i=1}^m\int_{[0,1]}\cdots\int_{[0,1]}\bigl (\sigma
_{r_{i}}(\cdot)
f-\sigma_{m+1}(\cdot) f \bigr)\,\d\mu_{r_{1}}\cdots\,\d\mu_{r_{m}}\,\d\mu
_{r} \Biggr] ,
\nonumber
\end{eqnarray}
which equals (\ref{gen2}) for appropriate values of $q,c,d,s$ and for
univariate $\sigma$.
The statement with $C_{\mathscr{P}(\mathbb{X})^{\#\mathcal
{I}}}([0,\infty))$ replaced by $D_{\mathscr{P}(\mathbb{X})^{\#
\mathcal{I}}}([0,\infty))$ now follows from
Theorems 1.6.1 and 4.2.11 of~\cite{EK86}, which, respectively, imply
the strong convergence of the corresponding semigroups and the weak
convergence of the law of $Y^{(n)}(\cdot)$ to that of $Y(\cdot)$.
Replacing $D_{\mathscr{P}(\mathbb{X})^{\#\mathcal{I}}}([0,\infty))$
with $C_{\mathscr{P}(\mathbb{X})^{\#\mathcal{I}}}([0,\infty))$
follows from~\cite{B68}, Section 18, by relativization of the Skorohod
topology to $C_{\mathscr{P}(\mathbb{X})^{\#\mathcal{I}}}([0,\infty))$.
\end{pf*}
\end{appendix}

\section*{Acknowledgements}
The authors are grateful to the Editor, an Associate Editor and a
referee for valuable remarks and suggestions that have lead to a
substantial improvement in the presentation. Thanks are also due to
Tommaso Frattini and Filippo Taddei for useful discussions.
This research was supported by the European Research Council (ERC)
through StG ``N-BNP'' 306406.

% imsref loaded by smiklovaite, 2011-12-19 07:13:47
%
% imsref loaded by smiklovaite, 2011-12-20 16:13:49

\printhistory

\end{document}